\input ./style/arxiv-general.cfg
\documentclass[aos,MSNbibl,nameyear,seceqn,dvips]{arximspdf}
\makeatletter
   \@ifpackageloaded{graphicx}{}{\usepackage{graphicx}}
\makeatother

\doi{10.1214/15-AOS1377}
\volume{44}
\issue{2}
\pubyear{2016}
\firstpage{564}
\lastpage{597}
\docsubty{FLA}

\makeatletter
\newcommand{\rrvert}{\vert}
\newcommand{\llvert}{\vert}
\newtheorem{theorem}{Theorem}[section]
\newtheorem{corollary}[theorem]{Corollary}
\newproclaim{example}{Example}
\newtheorem{lemma}[theorem]{Lemma}
\makeatother

\begin{document}
\begin{frontmatter}

\title{Optimal shrinkage estimation of mean parameters in family of
distributions with quadratic variance}
\runtitle{Optimal shrinkage estimation}

\begin{aug}
\author[A]{\fnms{Xianchao} \snm{Xie}\ead[label=e1]{xie1981@gmail.com}},
\author[B]{\fnms{S. C.} \snm{Kou}\corref{}\thanksref{T2}\ead[label=e2]{kou@stat.harvard.edu}}
\and
\author[C]{\fnms{Lawrence} \snm{Brown}\thanksref{T3}\ead[label=e3]{lbrown@wharton.upenn.edu}}
\runauthor{X. Xie, S. C. Kou and L. Brown}
\thankstext{T2}{Supported in part by grants from NIH and NSF.}
\thankstext{T3}{Supported in part by grants from NSF.}
\affiliation{Two Sigma Investments LLC, Harvard University and 
University of Pennsylvania}
\address[A]{X. Xie\\
Two Sigma Investments LLC\\
100 Avenue of the Americas, Floor 16\\
New York, New York 10013\\
USA\\
\printead{e1}}
\address[B]{S. C. Kou\\
Department of Statistics\\
Harvard University\\
Cambridge, Massachusetts 02138\\
USA\\
\printead{e2}}
\address[C]{L. Brown\\
The Wharton School\\
University of Pennsylvania\\
Philadelphia, Pennsylvania 19104\\
USA\\
\printead{e3}}
\end{aug}

%
\received{\smonth{2} \syear{2015}}
%
\revised{\smonth{8} \syear{2015}}

%
\begin{abstract}
This paper discusses the
simultaneous inference of mean parameters in a family of distributions
with quadratic variance function. We first introduce a class of
semiparametric/parametric shrinkage estimators and establish their
asymptotic optimality properties. Two specific cases, the location-scale
family and the natural exponential family with quadratic variance function,
are then studied in detail. We conduct a comprehensive simulation study to
compare the performance of the proposed methods with existing shrinkage
estimators. We also apply the method to real data and obtain encouraging
results.
\end{abstract}

%
%
\end{frontmatter}

\section{Introduction}\label{sec1}

The simultaneous inference of several mean parameters has interested
statisticians since the 1950s and has been widely applied in many scientific
and engineering problems ever since. \citet{Ste56} and \citet{JamSte61}
discussed the homoscedastic (equal variance) normal model and proved that
shrinkage estimators can have uniformly smaller risk compared to the
ordinary maximum likelihood estimate. This seminal work inspired a broad
interest in the study of shrinkage estimators in hierarchical normal models.
Efron and Morris (\citeyear{EfrMor72}, \citeyear{EfrMor73}) studied the James--Stein estimators in an
empirical Bayes framework and proposed several competing shrinkage
estimators. \citet{BerStr96} discussed this problem from a
hierarchical Bayesian perspective. For applications of shrinkage
techniques in practice,
see \citet{EfrMo75}, \citet{Rub81}, \citet{Mor83}, Green and
Strawderman (\citeyear{GrSt85}),
\citet{Jon91} and \citet{Bro08}.

There has also been substantial discussion on simultaneous inference for
non-Gaussian cases. \citet{Bro66} studied the admissibility of invariance
estimators in general location families. \citet{Joh84} discussed
inference in the context of Poisson estimation. \citet{CleZid75}
showed that Stein's effect also exists in the Poisson case, while a more
general treatment of discrete exponential families is given in \citet{Hwa82}.
The application of non-Gaussian hierarchical models has also flourished.
\citet{MosWal64} used a hierarchical Bayesian model based
on the
negative-binomial distribution to study the authorship of the Federalist
Papers. \citet{BerChr79} analyzed the binomial data using a
mixture of dirichlet processes. Gelman and Hill (\citeyear{GelHil07}) gave a comprehensive
discussion of data analysis using hierarchical models. Recently, nonparametric
empirical Bayes methods have been proposed by \citet{BroGre09},
Jiang and Zhang (\citeyear{JiaZha09}, \citeyear{JiaZha10})
and \citet{KoeMiz14}.

In this article, we focus on the simultaneous inference of the mean
parameters for families of distributions with quadratic variance function.
These distributions include many common ones, such as the normal, Poisson,
binomial, negative-binomial and gamma distributions; they also include
location-scale families, such as the $t$, logistic, uniform, Laplace, Pareto
and extreme value distributions. We ask the question: among all the
estimators that estimate the mean parameters by shrinking the within-group
sample mean toward a central location, is there an optimal one, subject to
the intuitive constraint that more shrinkage is applied to observations with
larger variances (or smaller sample sizes)? We propose a class of
semiparametric/parametric shrinkage estimators and show that there is
indeed an asymptotically optimal shrinkage estimator; this estimator is
explicitly obtained by minimizing an unbiased estimate of the risk. We note
that similar types of estimators are found in \citet{XieKouBro12} in
the context of the heteroscedastic (unequal variance) hierarchical
normal model.
The treatment in this article, however, is far more general, as it
covers a
much wider range of distributions. We illustrate the performance of our
shrinkage estimators by a comprehensive simulation study on both exponential
and location-scale families. We apply our shrinkage estimators on the
baseball data obtained by \citet{Bro08}, observing quite encouraging results.

The remainder of the article is organized as follows: In Section~\ref
{sec2}, we
introduce the general class of semiparametric URE shrinkage estimators,
identify the asymptotically optimal one and discuss its properties. We then
study the special case of location-scale families in Section~\ref{sec3}
and the case
of natural exponential families with quadratic variance (NEF-QVF) in
Section~\ref{sec4}. A systematic simulation study is conducted in
Section~\ref{sec5}, along with the
application of the proposed methods to the baseball data set in Section~\ref{sec6}. We
give a brief discussion in Section~\ref{sec7} and the technical proofs
are placed in the
\hyperref[app]{Appendix}.

\section{Semiparametric estimation of mean parameters of distributions with
quadratic variance function}
\label{sec2}

Consider simultaneously estimating the mean parameters of $p$ independent
observations $Y_{i}$, $i=1,\ldots,p$. Assume that the observation $Y_{i}$
comes from a distribution with quadratic variance function, that is, $%
E(Y_{i})=\theta_{i}\in\Theta$ and $\operatorname{Var}(Y_{i})=V(\theta_{i})/\tau_{i}$
such that
\[
V(\theta_{i})=\nu_{0}+\nu_{1}
\theta_{i}+\nu_{2}\theta_{i}^{2}
\]
with $\nu_{k}$ $(k=0,1,2)$ being known constants. The set $\Theta$ of
allowable parameters is a subset of $\{\theta:V(\theta)\geq0\}$.
$\tau
_{i}$ is assumed to be known here and can be interpreted as the within-group
sample size (i.e., when $Y_{i}$ is the sample average of the $i$th group)
or as the (square root) inverse-scale of $Y_{i}$. It is worth emphasizing
that distributions with quadratic variance function include many common
ones, such as the normal, Poisson, binomial, negative-binomial and gamma
distributions as well as location-scale families. We introduce the general
theory on distributions with quadratic variance function in this
section and
will specifically treat the cases of location-scale family and exponential
family in the next two sections.

In a simultaneous inference problem, hierarchical models are often used to
achieve partial pooling of information among different groups. For example,
in the famous normal-normal hierarchical model $Y_{i}\stackrel{\mathrm
{ind.}}{\sim}%
N(\theta_{i},1/\tau_{i})$, one often puts a conjugate prior
distribution $%
\theta_{i}\stackrel{\mathrm{i.i.d.}}{\sim}N(\mu,\lambda)$ and uses
the posterior
mean
%
\begin{equation}
E(\theta_{i}|\mathbf{Y};\mu,\lambda)=\frac{\tau_{i}}{\tau
_{i}+1/\lambda}%
\cdot
Y_{i}+\frac{1/\lambda}{\tau_{i}+1/\lambda}\cdot\mu \label{normal}
\end{equation}
to estimate $\theta_{i}$. Similarly, if $Y_{i}$ represents the within-group
average of Poisson observations $\tau_iY_{i}\stackrel{\mathrm
{ind.}}{\sim} \operatorname{Poisson}(\tau
_{i}\theta_{i}$), then with a conjugate gamma prior
distribution $\theta_{i}\stackrel{\mathrm{i.i.d.}}{\sim}\Gamma
(\alpha,\lambda)$,
the posterior mean
%
\begin{equation}
E(\theta_{i}|\mathbf{Y};\alpha,\lambda)=\frac{\tau_{i}}{\tau
_{i}+1/\lambda}%
\cdot Y_{i}+\frac{1/\lambda}{\tau_{i}+1/\lambda}\cdot\alpha\lambda \label{poisson}
\end{equation}
is often used to estimate $\theta_{i}$. The hyper-parameters, $(\mu
,\lambda)$ or $(\alpha,\lambda)$ above, are usually first estimated from
the marginal distribution of $Y_{i}$ and then plugged into the above
formulae to form an empirical Bayes estimate.

One potential drawback of the formal parametric empirical Bayes approaches
lies in its explicit parametric assumption on the prior distribution.
It can
lead to undesirable results if the explicit parametric assumption is
violated in real applications---we will see a real-data example in
Section~\ref{sec6}. Aiming to provide more flexible and, at the same
time, efficient
simultaneous estimation procedures, we consider in this section a class of
semiparametric shrinkage estimators.

To motivate these estimators, let us go back to the normal and Poisson
examples (\ref{normal}) and (\ref{poisson}). It is seen that the Bayes
estimate of each mean parameter $\theta_{i}$ is the weighted average
of $%
Y_{i}$ and the prior mean $\mu$ (or $\alpha\lambda$). In other
words, $%
\theta_{i}$ is estimated by shrinking $Y_{i}$ toward a central
location ($%
\mu$ or $\alpha\lambda$). It is also noteworthy that the amount of
shrinkage is governed by $\tau_{i}$, the sample size: the larger the sample
size, the less is the shrinkage toward the central location. This feature
makes intuitive sense. We will see in Section~\ref{sec4.2} that in fact these
observations hold not only for normal and Poisson distributions, but also
for general natural exponential families.

With these observations in mind, we consider in this section shrinkage
estimators of the following form:
%
\begin{equation}
\hat{\theta}{}^{\mathbf{b},\mu}_{i}=(1-b_{i})
\cdot{Y}_{i}+b_{i}\cdot\mu \label{semiform}
\end{equation}
with $b_{i}\in{}[0,1]$ satisfying
%
\begin{equation}
\mbox{Requirement (MON)}:b_{i}\leq b_{j}\mbox{ for any }i
\mbox{ and }j\mbox{ such that }\tau_{i}\geq\tau_{j}.
\label{Mon}
\end{equation}
Requirement (MON) asks the estimator to shrink the group mean with a larger
sample size (or smaller variance) less toward the central location $\mu$.
Other than this intuitive requirement, we do not put on any restriction
on $%
b_{i}$. Therefore, this class of estimators is semiparametic in nature.

The question we want to investigate is, for such a general estimator
$\hat{%
\bolds{\theta}}{}^{\mathbf{b},\mu}$, whether there exists an optimal
choice of $%
\mathbf{b}$ and $\mu$. Note that the two parametric estimates (\ref{normal})
and (\ref{poisson}) are special cases of the general class with
$b_{i}=\frac{%
1/\lambda}{\tau_{i}+1/\lambda}$. We will see shortly that such an optimal
choice indeed exists asymptotically (i.e., as $p\rightarrow\infty$) and
this asymptotically optimal choice is specified by an unbiased risk estimate
(URE).

For a general estimator $\hat{\bolds{\theta}}{}^{\mathbf{b},\mu}$ with
\emph{fixed} $%
\mathbf{b}$ and $\mu$, under the sum of squared-error loss
\[
l_{p}\bigl(\bolds{\theta},\hat{\bolds{\theta}}{}^{\mathbf{b},\mu}\bigr)=
\frac{1}{p} 
\sum_{i=1}^{p}
\bigl(\hat{\theta}{}^{\mathbf{b},\mu}_{i}-\theta_{i}
\bigr)^{2},
\]
an unbiased estimate of its risk $R_{p}(\bolds{\theta},\hat{\bolds
{\theta}}{}^{\mathbf{b}%
,\mu})=E[l_{p}(\bolds{\theta},\hat{\bolds{\theta}}{}^{\mathbf{b},\mu
})]$ is given by
\[
\operatorname{URE}(\mathbf{b},\mu)=\frac{1}{p}\sum_{i=1}^{p}
\biggl[ b_{i}^{2}\cdot(Y_{i}-\mu)^{2}+
( 1-2b_{i} ) \cdot\frac
{V(Y_{i})}{%
\tau_{i}+\nu_{2}} \biggr],
\]
because
\begin{eqnarray*}
E \bigl[ \operatorname{URE}(\mathbf{b},\mu) \bigr] &=&\frac{1}{p}\sum
_{i=1}^{p} \bigl\{ b_{i}^{2}
\bigl[\operatorname{Var}%
(Y_{i})+(\theta_{i}-
\mu)^{2}\bigr]+(1-2b_{i})\operatorname{Var}(Y_{i}) \bigr
\}
\\
&=&\frac{1}{p}\sum_{i=1}^{p}
\bigl[(1-b_{i})^{2}\operatorname{Var}(Y_{i})+b_{i}^{2}(
\theta_{i}-\mu)^{2}\bigr]=R_{p}\bigl(
\bolds{\theta},\hat{\bolds{\theta }}{}^{\mathbf{b},\mu}\bigr).
\end{eqnarray*}
Note that the idea and results below can be easily extended to the case of
weighted quadratic loss, with the only difference being that the regularity
conditions will then involve the corresponding weight sequence.

Ideally the ``best'' choice of $\mathbf{b}$
and $\mu$ is the one that minimizes $R_{p}(\bolds{\theta},\hat{\bolds
{\theta}}{}^{%
\mathbf{b},\mu})$, which is, however, unobtainable as the risk depends
on the
unknown $\bolds{\theta}$. To bypass this impracticability, we minimize
URE, the
unbiased estimate, with respect to $(\mathbf{b},\mu)$ instead. This
gives our
semiparametric URE shrinkage estimator:
%
\begin{equation}
\hat{\theta}{}^{\mathrm{SM}}_{i}=(1-\hat{b}_{i})\cdot
Y_{i}+\hat{b}_{i}\cdot\hat {\mu}%
{}^{\mathrm{SM}},
\label{semiURE}
\end{equation}
where
\begin{eqnarray*}
\bigl(\hat{\mathbf{b}}{}^{\mathrm{SM}},\hat{\mu}{}^{\mathrm{SM}}\bigr) &=&
\mbox{minimizer of }\operatorname {URE}(\mathbf{b}%
,\mu)
\\
&&\mbox{subject to }b_{i}\in{}[0,1], \mu\in{}\Bigl[ -\max
_{i}|Y_{i}|,\max_{i}|Y_{i}|
\Bigr]\cap\Theta
\\
&& \mbox{and Requirement (MON).}
\end{eqnarray*}
We require $\llvert \mu\rrvert \leq\max_{i}\llvert
Y_{i}\rrvert $, since no sensible estimators shrink the
observations to a location completely outside the range of the data.
Intuitively, the URE shrinkage estimator would behave well if $\operatorname
{URE}(%
\mathbf{b},\mu)$ is close to the risk $R_{p}(\bolds{\theta},\hat{\bolds
{\theta}}{}^{%
\mathbf{b},\mu})$.

To investigate the properties of the semiparametric URE shrinkage
estimator, we now introduce the following regularity conditions:
\begin{longlist}[(A)]
\item[(A)] $\limsup_{p\rightarrow\infty}\frac{1}{p}\sum_{i=1}^{p}%
\operatorname{Var}(Y_{i})<\infty$;

\item[(B)] $\limsup_{p\rightarrow\infty}\frac{1}{p}\sum_{i=1}^{p}%
\operatorname{Var}(Y_{i})\cdot\theta_{i}^{2}<\infty$;

\item[(C)] $\limsup_{p\rightarrow\infty}\frac{1}{p}\sum_{i=1}^{p}%
\operatorname{Var}(Y_{i}^{2})<\infty$;

\item[(D)] $\sup_{i}(\frac{\tau_{i}}{\tau_{i}+\nu_{2}})^{2}<\infty
$, $\sup_{i}(\frac{\nu_{1}}{\tau_{i}+\nu_{2}})^{2}<\infty$;

\item[(E)] $\limsup_{p\rightarrow\infty}\frac{1}{p^{1-\varepsilon}}%
E(\max_{1\leq i\leq p}Y_{i}^{2})<\infty$ for some $\varepsilon>0$.
\end{longlist}

The theorem below shows that $\operatorname{URE}(\mathbf{b},\mu)$ not only unbiasedly
estimates the risk, but also serves as a good approximation of the actual
loss $l_{p}(\bolds{\theta},\hat{\bolds{\theta}}{}^{\mathbf{b},\mu})$,
which is a much
stronger property. In fact, $\operatorname{URE}(\mathbf{b},\mu)$ is
asymptotically uniformly
close to the actual loss. Therefore, we expect that minimizing $\operatorname
{URE}(\mathbf{b},\mu)$
would lead to an estimate with competitive risk properties.

\begin{theorem}
\label{theorem1}Assuming regularity conditions \textup{(A)--(E)}, we have
\[
\sup\bigl\llvert \operatorname{URE}(\mathbf{b},\mu)-l_{p}\bigl(\bolds{\theta
},\hat{\bolds{\theta}}%
{}^{\mathbf{b},\mu}\bigr)\bigr\rrvert \rightarrow0\qquad
\mbox{in }L^{1}\mbox{ and in probability, as }p\rightarrow\infty,
\]
where the supremum is taken over $b_{i}\in{}[0,1]$, $\llvert
\mu
\rrvert \leq\max_{i}\llvert  Y_{i}\rrvert $ and Requirement (MON).
\end{theorem}

The following result compares the asymptotic behavior of our URE shrinkage
estimator (\ref{semiURE}) with other shrinkage estimators from the general
class. It establishes the asymptotic optimality of our URE shrinkage
estimator.

\begin{theorem}
\label{theorem2}Assume regularity conditions \textup{(A)--(E)}. Then for any shrinkage
estimator $\hat{\bolds{\theta}}{}^{\hat{\mathbf{b}},\hat{\mu}}=(1-\hat
{\mathbf{b}})\cdot{%
\mathbf{Y}}+\hat{\mathbf{b}}\cdot\hat{\mu}$, where $\hat{\mathbf
{b}}\in{}[
0,1] $ satisfies Requirement (MON), and $\llvert \hat{\mu}\rrvert
\leq\max_{i}\llvert  Y_{i}\rrvert $, we always have
\[
\lim_{p\rightarrow\infty}P \bigl( l_{p}\bigl(\bolds{\theta},\bolds{
\hat {\theta}}%
{}^{\mathrm{SM}}\bigr)\geq l_{p}\bigl(\bolds{
\theta},\hat{\bolds{\theta}}{}^{\hat{\mathbf
{b}},\hat{\mu}%
}\bigr)+\varepsilon \bigr) =0\qquad\mbox{for
any } \varepsilon>0
\]
and
\[
\limsup_{p\rightarrow\infty} \bigl[ R\bigl(\bolds{\theta},\hat{\bolds {
\theta}}%
{}^{\mathrm{SM}}\bigr)-R\bigl(\bolds{\theta},\hat{\bolds{
\theta}}{}^{\hat{\mathbf{b}},\hat{\mu
}}\bigr) \bigr] \leq0.
\]
\end{theorem}

As a special case of Theorem~\ref{theorem2}, the semiparametric URE
shrinkage estimator asymptotically dominates the parametric empirical Bayes
estimators, like (\ref{normal}) or (\ref{poisson}). It is worth noting that
the asymptotic optimality of our semiparametric URE shrinkage estimators
does not assume any prior distribution on the mean parameters $\theta_{i}$,
nor does it assume any parametric form on the distribution of $\mathbf{Y}$
(other than the quadratic variance function). Therefore, the results
enjoy a
large extent of robustness. In fact, the individual $Y_{i}$'s do not even
have to be from the same distribution family as long as the regularity
conditions (A)--(E) are met. (However, whether shrinkage estimation in that
case is a good idea or not becomes debatable.)

\subsection{Shrinking toward the grand mean}\label{sec2.1}

In the previous development, the central shrinkage location $\mu$ is
determined by minimizing URE. The joint minimization of $\mathbf{b}$
and $\mu$
offers asymptotic optimality in the class of estimators. For small or
moderate $p$ (the number of $Y_{i}$'s), however, it is not necessarily true
that the semiparametric URE shrinkage estimator will always be the optimal
one. In this setting, it might be beneficial to set $\mu$ by a
predetermined rule and only optimize~$\mathbf{b}$, as it might reduce the
variability of the resulting estimate. In this subsection, we consider
shrinking toward the grand mean:
\[
\hat{\mu}=\bar{Y}=\frac{1}{p}\sum_{i=1}^{p}{Y}_{i}.
\]

The particular reason why the grand average is chosen instead of the
weighted average
$\bar{Y}_w=(\sum_{i=1}^{p}\tau_i Y_{i})/(\sum_{i=1}^{p}\tau_{i})$
is that the latter might be biased when $\theta_{i}$ and $\tau_{i}$
are dependent. In the case where such dependence is not a concern, the
idea and results
obtained below can be similarly derived.

The corresponding estimator becomes
%
\begin{equation}
\hat{\theta}{}^{\mathbf{b},\bar{Y}}_{i}=(1-b_{i})Y_{i}+b_{i}{
\bar{Y}}, \label{UREgrandmean}
\end{equation}
where $b_{i}\in{}[0,1]$ satisfies Requirement (MON). To find the
asymptotically optimal choice of $\mathbf{b}$, we start from an unbiased
estimate of the risk of $\hat{\bolds{\theta}}{}^{\mathbf{b},\bar{Y}}$.

It is straightforward to verify that for fixed $\mathbf{b}$ an unbiased estimate
of the risk of $\hat{\bolds{\theta}}{}^{\mathbf{b},{\bar{Y}}}$ is
\[
\operatorname{URE}^{G}(\mathbf{b})=\frac{1}{p}\sum
_{i=1}^{p} \biggl[ b_{i}^{2}(Y_{i}-
\bar{Y})^{2}+ \biggl( 1-2\biggl(1-\frac{1}{p}\biggr)b_{i}
\biggr) \frac
{%
V(Y_{i})}{\tau_{i}+\nu_{2}} \biggr].
\]
Note that we use the superscript ``$G$'',
which stands for ``grand mean'', to
distinguish it from the previous $\operatorname{URE}(\mathbf{b},\mu)$. Like
what we
did previously, minimizing the URE$^{G}$ with respect to $\mathbf{b}$
then leads
to our semiparametric URE grand-mean shrinkage estimator
%
\begin{equation}
\hat{\theta}{}^{\mathrm{SG}}_{i}=\bigl(1-\hat{b}{}^{\mathrm{SG}}_{i}
\bigr)\cdot Y_{i}+\hat {b}{}^{\mathrm{SG}}_{i}\cdot \bar{Y},
\label{grandmeanspshrink}
\end{equation}
where
\begin{eqnarray*}
\hat{\mathbf{b}}{}^{\mathrm{SG}} &=&\mbox{minimizer of }\operatorname{URE}^{G}(
\mathbf {b})
\\
&&\mbox{subject to }b_{i}\in{}[0,1]\mbox{ and Requirement (MON).}
\end{eqnarray*}

Again, we expect that the URE estimator $\hat{\bolds{\theta}}{}^{\mathrm{SG}}$
would be
competitive if $\operatorname{URE}^{G}$ is close to the risk function or the loss
function. The next theorem confirms the uniform closeness.

\begin{theorem}
\label{theorem3} Under regularity conditions \textup{(A)--(E)}, we have
\[
\sup\bigl\llvert \operatorname{URE}^{G}(\mathbf{b})-l_{p}\bigl(
\bolds{\theta},\bolds {\hat{\theta}}%
{}^{\mathbf{b},{\bar{Y}}}\bigr)\bigr\rrvert
\rightarrow0\qquad \mbox{in }L^{1}\mbox{ and in probability, as }p\rightarrow
\infty,
\]
where the supremum is taken over $b_{i}\in{}[0,1]$ and Requirement
(MON).
\end{theorem}

Consequently, $\hat{\bolds{\theta}}{}^{\mathrm{SG}}$ is asymptotically optimal
among all
shrinkage estimators $\hat{\bolds{\theta}}{}^{\mathbf{b},\bar{Y}}$ that
shrink toward
the grand mean ${\bar{Y}}$, as shown in the next theorem.

\begin{theorem}
\label{theorem4}Assume regularity conditions \textup{(A)--(E)}. Then for any shrinkage
estimator $\hat{\bolds{\theta}}{}^{\hat{\mathbf{b}},{\bar{Y}}}=(1-\hat
{\mathbf{b}})\cdot
\mathbf{Y}+\hat{\mathbf{b}}\cdot{\bar{Y}}$, where $\hat{\mathbf{b}}\in
{}[0,1]$
satisfies Requirement (MON), we have
\[
\lim_{p\rightarrow\infty}P \bigl( l_{p}\bigl(\bolds{\theta},\bolds{
\hat {\theta}}%
{}^{\mathrm{SG}}\bigr)\geq l_{p}\bigl(\bolds{
\theta},\hat{\bolds{\theta}}{}^{\hat{\mathbf
{b}},\bar{Y}%
}\bigr)+\varepsilon \bigr) =0\qquad \mbox{for
any } \varepsilon>0
\]
and
\[
\limsup_{p\rightarrow\infty} \bigl[ R\bigl(\bolds{\theta},\hat{\bolds {
\theta}}%
{}^{\mathrm{SG}}\bigr)-R\bigl(\bolds{\theta},\hat{\bolds{
\theta}}{}^{\hat{\mathbf{b}},\bar
{Y}}\bigr) \bigr] \leq 0.
\]
\end{theorem}

\section{Simultaneous estimation of mean parameters in location-scale
families}\label{sec3}

In this section, we focus on location-scale families, which are a special
case of distributions with quadratic variance functions. We show how the
regularity conditions can be simplified in this case.

For a location-scale family, we can write
%
\begin{equation}
Y_{i}=\theta_{i}+Z_{i}/\sqrt{
\tau_{i}}, \label{loc-scale}
\end{equation}
where the standard variates $Z_{i}$ are i.i.d. with mean zero and
variance $%
\nu_{0}$. The constants $\nu_{1}$ and $\nu_{2}$ in the quadratic variance
function $V(\theta_{i})=\nu_{0}+\nu_{1}\theta_{i}+\nu_{2}\theta_{i}^{2}
$ are zero for the location-scale family. $1/\sqrt{\tau_{i}}$ is the scale
of $Y_{i}$.

The next lemma simplifies the regularity conditions for a location-scale
family.

\begin{lemma}
\label{theoremLS} For $Y_{i}$, $i=1,\ldots,p$, independently from a
location-scale family (\ref{loc-scale}), the following four conditions imply
the regularity conditions \textup{(A)--(E)} in Section~\ref{sec2}:
\begin{enumerate}[(iii)]
\item[(i)] $\limsup_{p\rightarrow\infty}\frac{1}{p}%
\sum_{i=1}^{p}1/\tau_{i}^{2}<\infty$;

\item[(ii)] $\limsup_{p\rightarrow\infty}\frac{1}{p}%
\sum_{i=1}^{p}\theta_i^{2}/\tau_{i}<\infty$;

\item[(iii)] $\limsup_{p\rightarrow\infty}\frac{1}{p}%
\sum_{i=1}^{p}\llvert {\theta_i}\rrvert ^{2+\varepsilon}<\infty$
for some $\varepsilon>0$;

\item[(iv)] the standard variate $Z$ satisfies $P(\llvert  Z\rrvert
>t)\leq Dt^{-\alpha}$ for constants $D>0$, $\alpha>4$.
\end{enumerate}
\end{lemma}

Note that (i)--(iv) in Lemma~\ref{theoremLS} covers the common location-scale
families, including the $t$ (degree of freedom $>4$), normal, uniform,
logistic, Laplace, Pareto ($\alpha>4$) and extreme value distributions.

Lemma~\ref{theoremLS}, together with Theorems \ref{theorem1} and \ref%
{theorem2}, immediately yields the following corollaries: the
semiparametric URE shrinkage estimator is asymptotically optimal for
location-scale families.

\begin{corollary}
For $Y_{i}$, $i=1,\ldots,p$, independently from a location-scale
family (%
\ref{loc-scale}), under conditions \textup{(i)--(iv)} in Lemma~\ref{theoremLS}, we
have
\[
\sup\bigl\llvert \operatorname{URE}(\mathbf{b},\mu)-l_{p}\bigl(\bolds{\theta
},\hat{\bolds{\theta}}%
{}^{\mathbf{b},\mu}\bigr)\bigr\rrvert \rightarrow0\qquad
\mbox{in }L^{1}\mbox{ and in probability, as }p\rightarrow\infty,
\]
where the supremum is taken over $b_{i}\in{}[0,1]$, $\llvert
\mu
\rrvert \leq\max_{i}\llvert  Y_{i}\rrvert $ and Requirement (MON).
\end{corollary}

\begin{corollary}
Let $Y_{i}$, $i=1,\ldots,p$, be independent from a location-scale
family (%
\ref{loc-scale}). Assume conditions \textup{(i)--(iv)} in Lemma~\ref{theoremLS}. Then
for any shrinkage estimator $\hat{\bolds{\theta}}{}^{\hat{\mathbf{b}},\hat
{\mu}}=(1-%
\hat{\mathbf{b}})\cdot{\mathbf{Y}}+\hat{\mathbf{b}}\cdot\hat{\mu}$,
where $\hat{%
\mathbf{b}}\in{}[0,1]$ satisfies Requirement (MON), and $\llvert \hat{%
\mu}\rrvert \leq\max_{i}\llvert  Y_{i}\rrvert $, we always have
\[
\lim_{p\rightarrow\infty}P \bigl( l_{p}\bigl(\bolds{\theta},\bolds{
\hat {\theta}}%
{}^{\mathrm{SM}}\bigr)\geq l_{p}\bigl(\bolds{
\theta},\hat{\bolds{\theta}}{}^{\hat{\mathbf
{b}},\hat{\mu}%
}\bigr)+\varepsilon \bigr) =0\qquad\mbox{for
any } \varepsilon>0
\]
and
\[
\limsup_{p\rightarrow\infty} \bigl[ R\bigl(\bolds{\theta},\hat{\bolds {
\theta}}%
{}^{\mathrm{SM}}\bigr)-R\bigl(\bolds{\theta},\hat{\bolds{
\theta}}{}^{\hat{\mathbf{b}},\hat{\mu
}}\bigr) \bigr] \leq0.
\]
\end{corollary}

In the case of shrinking toward the grand mean ${\bar{Y}}$, the
corresponding semiparametric URE grand-mean shrinkage estimator is also
asymptotically optimal.

\begin{corollary}
For $Y_{i}$, $i=1,\ldots,p$, independently from a location-scale
family (%
\ref{loc-scale}), under conditions \textup{(i)--(iv)} in Lemma~\ref{theoremLS}, we
have
\[
\sup\bigl\llvert \operatorname{URE}^{G}(\mathbf{b})-l_{p}\bigl(
\bolds{\theta},\bolds {\hat{\theta}}%
{}^{\mathbf{b},{\bar{Y}}}\bigr)\bigr\rrvert
\rightarrow0\qquad\mbox{in }L^{1}\mbox{ and in probability, as }p\rightarrow
\infty,
\]
where the supremum is taken over $b_{i}\in{}[0,1]$ and Requirement
(MON).
\end{corollary}

\begin{corollary}
Let $Y_{i}$, $i=1,\ldots,p$, be independent from a location-scale
family (%
\ref{loc-scale}). Assume conditions \textup{(i)--(iv)} in Lemma~\ref{theoremLS}. Then
for any shrinkage estimator $\hat{\bolds{\theta}}{}^{\hat{\mathbf
{b}},{\bar{Y}}}=(1-%
\hat{\mathbf{b}})\cdot\mathbf{Y}+\hat{\mathbf{b}}\cdot{\bar{Y}}$,
where $\hat{%
\mathbf{b}}\in{}[0,1]$ satisfies Requirement (MON), we have
\[
\lim_{p\rightarrow\infty}P \bigl( l_{p}\bigl(\bolds{\theta},\bolds{
\hat {\theta}}%
{}^{\mathrm{SG}}\bigr)\geq l_{p}\bigl(\bolds{
\theta},\hat{\bolds{\theta}}{}^{\hat{\mathbf
{b}},\bar{Y}%
}\bigr)+\varepsilon \bigr) =0\qquad\mbox{for
any } \varepsilon>0
\]
and
\[
\limsup_{p\rightarrow\infty} \bigl[ R\bigl(\bolds{\theta},\hat{\bolds {
\theta}}%
{}^{\mathrm{SG}}\bigr)-R\bigl(\bolds{\theta},\hat{\bolds{
\theta}}{}^{\hat{\mathbf{b}},\bar
{Y}}\bigr) \bigr] \leq 0.
\]
\end{corollary}

\section{Simultaneous estimation of mean parameters in natural exponential
family with quadratic variance function}\label{sec4}

\subsection{Semiparametric URE shrinkage estimators}\label{sec4.1}

In this section, we focus on natural exponential families with quadratic
variance functions (NEF-QVF), as they incorporate the most common
distributions that one encounters in practice. We show how the regularity
conditions (A)--(E) can be significantly simplified and offer concrete
examples.

It is well known that there are in total six distinct distributions that
belong to NEF-QVF [\citet{Mor82}]: the normal, binomial, Poisson,
negative-binomial, Gamma and generalized hyperbolic secant (GHS)
distributions. We represent in general an NEF-QVF as
\[
Y_{i}\sim\mbox{NEF-QVF}\bigl[\theta_{i},V(
\theta_{i})/\tau_{i}\bigr],
\]
where $\theta_{i}=E(Y_{i})\in\Theta$ is the mean parameter and $\tau_{i}$
is the convolution parameter (or within-group sample size). For
example, in
the binomial case $Y_{i}\sim \operatorname{Bin}(n_{i},p_{i})/n_{i}$, $\theta
_{i}=p_{i}$, $%
V(\theta_{i})=\theta_{i}(1-\theta_{i})$ and $\tau_{i}=n_{i}$.

The next result provides easy-to-check conditions that considerably
simplify those in Section~\ref{sec2}. As the case of heteroscedastic normal
data has been studied in \citet{XieKouBro12}, we concentrate on the
other five NEF-QVF distributions.

\begin{lemma}
\label{lemmaNEFQVF} For the five non-Gaussian NEF-QVF distributions,
Table~\ref{tab1} lists the respective conditions, under which regularity
conditions \textup{(A)--(E)} in Section~\ref{sec2} are satisfied. For example,
for the
binomial distribution, the condition is $\tau_{i}=n_{i}\geq2$ for all $i$.
\end{lemma}

\begin{table}
\tabcolsep=0pt
\caption{The conditions for the five nonnormal NEF-QVF distributions to
guarantee the regularity conditions \textup{(A)--(E)} in Section \protect\ref{sec2}}
\label{tab1}
\begin{tabular*}{\textwidth}{@{\extracolsep{\fill}}lcccl@{}}
\hline
\textbf{Distribution} & \textbf{Data} & $\bolds{(\nu_{0},\nu_{1},\nu_{2})}$ & \textbf{Note} &
\multicolumn{1}{c@{}}{\textbf{Conditions}} \\
\hline
Binomial & $Y_{i}\sim \operatorname{Bin}(n_{i},p_{i})/n_{i}$ & $(0,1,-1)$ & $\tau
_{i}=n_{i} $ & $n_{i}\geq2$ for all $i=1,\ldots,p$
\\ [3pt]
Poisson & $Y_{i}\sim \operatorname{Poi}(\tau_{i}\theta_{i})/\tau_{i}$ & $(0,1,0)$ & &
(i) $\inf_{i}\tau_{i}>0$, $\inf_{i}(\tau
_{i}\theta
_{i})>0$ \\
& & & & (ii) $\sum_{i}\theta_{i}^{3}=O(p)$ \\[3pt]
Neg-binomial & $Y_{i}\sim N\operatorname{Bin}(n_{i},p_{i})/n_{i}$ &
 $(0,1,1)$ & $\tau_{i}=n_{i}$ &
 (i) $\inf_{i}(n_{i}p_{i})>0$ \\
& & & $\theta_{i}=\frac{p_{i}}{1-p_{i}}$ &
(ii) $\sum_{i}(\frac{p_{i}}{1-p_{i}})^{4}=O(p)$ \\ [3pt]
Gamma & $Y_{i}\sim\Gamma(\tau_{i}\alpha,\lambda_{i})/\tau_{i}$ & $
(0,0,1/\alpha)$ & $\theta_{i}=\alpha\lambda_{i}$ &
(i) $\inf_{i}\tau_{i}>0$ \\
& & & & (ii) $\sum_{i}\lambda_{i}^{4}=O(p)$ \\[3pt]
GHS & $Y_{i}\sim \operatorname{GHS}(\tau_{i}\alpha,\lambda_{i})/\tau_{i}$ &$(\alpha
,0,1/\alpha)$ & $\theta_{i}=\alpha\lambda_{i}$ &
(i) $\inf_{i}\tau_{i}>0$ \\
& & & & (ii) $\sum_{i}\lambda_{i}^{4}=O(p)$\\
\hline
\end{tabular*}
\end{table}

Lemma~\ref{lemmaNEFQVF} and the general theory in Section~\ref{sec2}
yield the following optimality results for our semiparametric URE shrinkage
estimator in the case of NEF-QVF.

\begin{corollary}
Let $Y_{i}\stackrel{{ind.}}{\sim}\mathrm{NEF\mbox{-}QVF}[\theta
_{i},V(\theta_{i})/\tau
_{i}]$, $i=1,\ldots,p$, be non-Gaussian. Under the respective conditions
listed in Table~\ref{tab1}, we have
\[
\sup\bigl\llvert \operatorname{URE}(\mathbf{b},\mu)-l_{p}\bigl(\bolds{\theta
},\hat{\bolds{\theta}}%
{}^{\mathbf{b},\mu}\bigr)\bigr\rrvert \rightarrow0\qquad
\mbox{in }L^{1}\mbox{ and in probability, as }p\rightarrow\infty,
\]
where the supremum is taken over $b_{i}\in{}[0,1]$, $\llvert
\mu
\rrvert \leq\max_{i}\llvert  Y_{i}\rrvert $ and Requirement (MON).
\end{corollary}

\begin{corollary}
Let $Y_{i}\stackrel{{ind.}}{\sim}\mathrm{NEF\mbox{-}QVF}[\theta
_{i},V(\theta_{i})/\tau
_{i}]$, $i=1,\ldots,p$, be non-Gaussian. Assume the respective conditions
listed in Table~\ref{tab1}. Then for any shrinkage estimator $%
\hat{\bolds{\theta}}{}^{\hat{\mathbf{b}},\hat{\mu}}=(1-\hat{\mathbf
{b}})\cdot{\mathbf{Y}}%
+\hat{\mathbf{b}}\cdot\hat{\mu}$, where $\hat{\mathbf{b}}\in{}[0,1]$
satisfies Requirement (MON), and $\llvert \hat{\mu}\rrvert \leq
\max_{i}\llvert  Y_{i}\rrvert $, we always have
\[
\lim_{p\rightarrow\infty}P \bigl( l_{p}\bigl(\bolds{\theta},\bolds{
\hat {\theta}}%
{}^{\mathrm{SM}}\bigr)\geq l_{p}\bigl(\bolds{
\theta},\hat{\bolds{\theta}}{}^{\hat{\mathbf
{b}},\hat{\mu}%
}\bigr)+\varepsilon \bigr) =0\qquad\mbox{for
any } \varepsilon>0
\]
and
\[
\limsup_{p\rightarrow\infty} \bigl[ R\bigl(\bolds{\theta},\hat{\bolds {
\theta}}%
{}^{\mathrm{SM}}\bigr)-R\bigl(\bolds{\theta},\hat{\bolds{
\theta}}{}^{\hat{\mathbf{b}},\hat{\mu
}}\bigr) \bigr] \leq0.
\]
\end{corollary}

For shrinking toward the grand mean ${\bar{Y}}$, the corresponding
semiparametric URE grand mean shrinkage estimator is also asymptotically
optimal (within the smaller class).

\begin{corollary}
Let $Y_{i}\stackrel{{ind.}}{\sim}\mathrm{NEF\mbox{-}QVF}[\theta
_{i},V(\theta_{i})/\tau
_{i}]$, $i=1,\ldots,p$, be non-Gaussian. Under the respective conditions
listed in Table~\ref{tab1}, we have
\[
\sup\bigl\llvert \operatorname{URE}^{G}(\mathbf{b})-l_{p}\bigl(
\bolds{\theta},\bolds {\hat{\theta}}%
{}^{\mathbf{b},{\bar{Y}}}\bigr)\bigr\rrvert
\rightarrow0\qquad\mbox{in }L^{1}\mbox{ and in probability, as }p\rightarrow
\infty,
\]
where the supremum is taken over $b_{i}\in{}[0,1]$ and Requirement
(MON).
\end{corollary}

\begin{corollary}
Let $Y_{i}\stackrel{{ind.}}{\sim}\mathrm{NEF\mbox{-}QVF}[\theta
_{i},V(\theta_{i})/\tau
_{i}]$, $i=1,\ldots,p$, be non-Gaussian. Assume the respective conditions
listed in Table~\ref{tab1}. Then for any shrinkage estimator $%
\hat{\bolds{\theta}}{}^{\hat{\mathbf{b}},{\bar{Y}}}=(1-\hat{\mathbf
{b}})\cdot\mathbf{Y}+%
\hat{\mathbf{b}}\cdot{\bar{Y}}$, where $\hat{\mathbf{b}}\in{}[
0,1]$ satisfies
Requirement (MON), we have
\[
\lim_{p\rightarrow\infty}P \bigl( l_{p}\bigl(\bolds{\theta},\bolds{
\hat {\theta}}%
{}^{\mathrm{SG}}\bigr)\geq l_{p}\bigl(\bolds{
\theta},\hat{\bolds{\theta}}{}^{\hat{\mathbf
{b}},\bar{Y}%
}\bigr)+\varepsilon \bigr) =0\qquad\mbox{for
any } \varepsilon>0
\]
and
\[
\limsup_{p\rightarrow\infty} \bigl[ R\bigl(\bolds{\theta},\hat{\bolds {
\theta}}%
^{\mathrm{SG}}\bigr)-R\bigl(\bolds{\theta},\hat{\bolds{
\theta}}{}^{\hat{\mathbf{b}},\bar
{Y}}\bigr) \bigr] \leq 0.
\]
\end{corollary}

\subsection{Parametric URE shrinkage estimators and conjugate priors}\label{sec4.2}

For $Y_{i}$ from an exponential family, hierarchical models based on the
conjugate prior distributions are often used; the hyper-parameters in the
prior distribution are often estimated from the marginal distribution
of $%
Y_{i}$ in an empirical Bayes way. Two questions arise naturally in this
scenario. First, are there other choices to estimate the hyper-parameters?
Second, is there an optimal choice? We will show in this subsection
that our
URE shrinkage idea applies to the parametric conjugate prior case; the
resulting \emph{parametric} URE shrinkage estimators are asymptotically
optimal, and thus asymptotically dominate the traditional empirical Bayes
estimators.

Let $Y_{i}\sim\mathrm{NEF\mbox{-}QVF}[\theta_{i},V(\theta_{i})/\tau_{i}]$,
$i=1,\ldots
,p$, be independent. If $\theta_{i}$ are i.i.d. from the conjugate prior,
the Bayesian estimate of $\theta_{i}$ is then
%
\begin{equation}
\hat{\theta}{}^{\gamma,\mu}_{i}=\frac{\tau_{i}}{\tau_{i}+\gamma}\cdot
Y_{i}+\frac{\gamma}{\tau_{i}+\gamma}\cdot\mu, \label{postmean}
\end{equation}
where $\gamma$ and $\mu$ are functions of the hyper-parameters in the prior
distribution. Table~\ref{tab2} details the conjugate priors for the five
well-known NEF-QVF distributions---the normal, binomial, Poisson,
negative-binomial and gamma distributions---and the corresponding
expressions of $\gamma$ and $\mu$ in terms of the hyper-parameters. For
example, in the binomial case, $Y_{i}\stackrel{\mathrm{ind.}}{\sim
}\operatorname{Bin}(\tau
_{i},\theta_{i})/\tau_{i}$ and the conjugate prior is $\theta
_{i}\stackrel{%
\mathrm{i.i.d.}}{\sim}\operatorname{Beta}(\alpha,\beta)$; $\gamma
=\alpha+\beta$, $\mu=\alpha
/(\alpha+\beta)$. Although it can be shown that for the sixth NEF-QVF
distribution---the GHS distribution---taking a conjugate prior also
gives (%
\ref{postmean}), the conjugate prior distribution does not have a clean
expression and is rarely encountered in practice. We thus omit it from
Table~\ref{tab2}.

\begin{table}
\tabcolsep=0pt
\caption{Conjugate priors for the five well-known NEF-QVF distributions}
\label{tab2}
\begin{tabular*}{\textwidth}{@{\extracolsep{\fill}}lccc@{}}
\hline
\textbf{Distribution} & \textbf{Data} & \textbf{Conjugate prior} & \multicolumn{1}{c@{}}{$\bolds{\gamma}$ \textbf{and} $\bolds{\mu}$} \\
\hline
Normal & $Y_{i}\sim N(\theta_{i},1/\tau_{i})$ & {$
\theta_{i}\stackrel{\mathrm{i.i.d.}}{\sim}N(\mu,\lambda)$} &
$\gamma=1/\lambda$
\\
Binomial & $Y_{i}\sim \operatorname{Bin}(\tau_{i},\theta_{i})/\tau_{i}$ &
{$\theta_{i}\stackrel{\mathrm{i.i.d.}}{\sim
}\operatorname{Beta}(\alpha,\beta)$%
} & $\gamma=\alpha+\beta$, $\mu=\frac{\alpha}{\alpha+\beta}$ \\
Poisson & $Y_{i}\sim \operatorname{Poi}(\tau_{i}\theta_{i})/\tau_{i}$ &
{$\theta_{i}\stackrel{\mathrm{i.i.d.}}{\sim
}\Gamma(\alpha
,\lambda)$} & $\gamma=1/\lambda$, $\mu=\alpha\lambda$ \\
Neg-binomial & $Y_{i}\sim N\operatorname{Bin}(\tau_{i},p_{i})/\tau_{i}$ &
{$p_{i}\stackrel{\mathrm{i.i.d.}}{\sim
}\operatorname{Beta}(\alpha,\beta)$, $%
\theta_{i}=\frac{p_{i}}{1-p_{i}}$} & $\gamma=\beta-1$, $\mu=\frac
{\alpha
}{\beta-1}$ \\
Gamma & $Y_{i}\sim\Gamma(\tau_{i}\alpha,\lambda_{i})/\tau_{i}$ &
{$\lambda_{i}\stackrel{\mathrm{i.i.d.}}{\sim
}\mathrm{inv}\mbox{-}\Gamma
(\alpha_{0},\beta_{0})$, $\theta_{i}=\alpha\lambda_{i}$} & $\gamma
=%
\frac{\alpha_{0}-1}{\alpha}$, $\mu=\frac{\alpha\beta_{0}}{\alpha_{0}-1}
$ \\
\hline
\end{tabular*}
\end{table}

We now apply our URE idea to formula (\ref{postmean}) to estimate
$(\gamma
,\mu)$, in contrast to the conventional empirical Bayes method that
determines the hyper-parameters through the marginal distribution. For
fixed $%
\gamma$ and $\mu$, an unbiased estimate for the risk of $\hat{\bolds
{\theta}}%
{}^{\gamma,\mu}$ is given by
\[
\operatorname{URE}^{P}(\gamma,\mu)=\frac{1}{p}\sum
_{i=1}^{p} \biggl[ \frac{%
\gamma^{2}}{(\tau_{i}+\gamma)^{2}}
\cdot(Y_{i}-\mu)^{2}+\frac{\tau
_{i}-\gamma}{\tau_{i}+\gamma}\cdot
\frac{V(Y_{i})}{\tau_{i}+\nu
_{2}}%
 \biggr] ,
\]
where we use the superscript ``$P$'' to
stand for ``parametric''. Minimizing $%
\operatorname{URE}^{P}(\gamma, \mu)$ leads to our parametric URE shrinkage
estimator
%
\begin{equation}
\hat{\theta}{}^{{\mathrm{PM}}}_{i}=\frac{\tau_{i}}{\tau_{i}+\hat{\gamma
}{}^{{\mathrm{PM}}}}\cdot
Y_{i}+\frac{\hat{\gamma}{}^{P{M}}}{\tau_{i}+\hat{\gamma}{}^{{\mathrm{PM}}}}\cdot \hat{\mu%
}{}^{{\mathrm{PM}}},
\label{parURE}
\end{equation}
where
\begin{eqnarray*}
\bigl(\hat{\gamma}{}^{{\mathrm{PM}}},\hat{\mu}{}^{\mathrm{PM}}\bigr) &=&\arg\min
\operatorname {URE}^{P}(\gamma ,\mu)
\\
&&\mbox{over }\Bigl\{0\leq\gamma\leq\infty,\llvert \mu\rrvert \leq \max
_{i}\llvert Y_{i}\rrvert ,\mu\in\Theta]\Bigr\}.
\end{eqnarray*}

Parallel to Theorems \ref{theorem1} and \ref{theorem2}, the next two results
show that the parametric URE shrinkage estimator gives the asymptotically
optimal choice of $(\gamma,\mu)$, if one wants to use estimators of the
form (\ref{postmean}).

\begin{theorem}
\label{theorem5}Let $Y_{i}\stackrel{{ind.}}{\sim}\mathrm
{NEF\mbox{-}QVF}[\theta
_{i},V(\theta_{i})/\tau_{i}]$, $i=1,\ldots,p$, be non-Gaussian.
Under the
respective conditions listed in Table~\ref{tab1}, we have
\[
\sup\bigl\llvert \operatorname{URE}^{P}(\gamma,\mu)-l_{p}\bigl(
\bolds{\theta},\hat {%
\bolds{\theta}}{}^{\gamma,\mu}\bigr)\bigr\rrvert
\rightarrow\qquad\mbox{in }L^{1}\mbox{ and in probability, as }p\rightarrow
\infty,
\]
where the supremum is taken over $\{0\leq\gamma\leq\infty,\llvert \mu
\rrvert \leq\max_{i}\llvert  Y_{i}\rrvert ,\mu\in
\Theta]\}$.
\end{theorem}

\begin{theorem}
\label{theorem6}Let $Y_{i}\stackrel{{ind.}}{\sim}\mathrm
{NEF\mbox{-}QVF}[\theta
_{i},V(\theta_{i})/\tau_{i}]$, $i=1,\ldots,p$, be non-Gaussian. Assume
the respective conditions listed in Table~\ref{tab1}. Then for any estimator
$\hat{\bolds{\theta}}{}^{\hat{\gamma},\hat{\mu}}=\frac{\bolds{\tau
}}{\bolds{\tau}+\hat{%
\gamma}}\mathbf{Y}+\frac{\hat{\gamma}}{\bolds{\tau}+\hat{\gamma}}\hat
{\mu}$, where $%
\hat{\gamma}\geq0$ and $\llvert \hat{\mu}\rrvert \leq
\max_{i}\llvert  Y_{i}\rrvert $, we always have
\[
\lim_{p\rightarrow\infty}P \bigl( l_{p}\bigl(\bolds{\theta},\hat{
\bolds {\theta}}%
{}^{\mathrm{PM}}\bigr)\geq l_{p}\bigl(\bolds{
\theta},\hat{\bolds{\theta}}{}^{\hat{\gamma},\hat
{\mu}%
}\bigr)+\varepsilon \bigr) =0\qquad\mbox{for
any } \varepsilon>0
\]
and
\[
\limsup_{p\rightarrow\infty} \bigl[ R_{p}\bigl(\bolds{\theta},
\hat{%
\bolds{\theta}}{}^{\mathrm{PM}}\bigr)-R_{p}\bigl(\bolds{
\theta},\hat{\bolds{\theta}}{}^{\hat
{\gamma},\hat{\mu%
}}\bigr) \bigr] \leq0.
\]
\end{theorem}

In the case of shrinking toward the grand mean ${\bar{Y}}$ (when $p$, the
number of $Y_{i}$'s, is small or moderate), we have the following
parametric results parallel to the semiparametric ones.

First, for the grand-mean shrinkage estimator
%
\begin{equation}
\hat{\theta}{}^{\gamma,\bar{Y}}=\frac{\tau_{i}}{\tau_{i}+\gamma}\cdot Y_{i}+
\frac{\gamma}{\tau_{i}+\gamma}\cdot\bar{Y}, \label{meanshrink}
\end{equation}
with a fixed $\gamma$, an unbiased estimate of its risk is
\[
\operatorname{URE}^{\mathrm{PG}}(\gamma)=\frac{1}{p}\sum
_{i=1}^{p} \biggl[ \frac{%
\gamma^{2}}{(\tau_{i}+\gamma)^{2}} (
Y_{i}-\bar{Y} ) ^{2}+ \biggl( 1-2 \biggl( 1-
\frac{1}{p} \biggr) \frac{\gamma}{\tau_{i}+\gamma} \biggr) \frac{V(Y_{i})}{\tau_{i}+\nu_{2}} \biggr].
\]

Minimizing it yields our parametric URE grand-mean shrinkage estimator
%
\begin{equation}
\hat{\theta}_{i}^{{\mathrm{PG}}}=\frac{\tau_{i}}{\tau_{i}+\hat{\gamma
}{}^{{\mathrm{PG}}}}\cdot
Y_{i}+\frac{\hat{\gamma}{}^{\mathrm{P{G}}}}{\tau_{i}+\hat{\gamma}{}^{{\mathrm{PG}}}}\cdot \bar{Y}, \label{meanUREshrink}
\end{equation}
where
\[
\hat{\gamma}{}^{{\mathrm{PG}}}=\arg\min_{0\leq\gamma\leq\infty}
\operatorname{URE}%
^{\mathrm{PG}}(\gamma).
\]

Similar to Theorems \ref{theorem5} and \ref{theorem6}, the next two results
show that in the case of shrinking toward the grand mean under the
formula (%
\ref{meanshrink}), the parametric URE grand-mean shrinkage estimator is
asymptotically optimal.

\begin{theorem}
\label{theorem7}Let $Y_{i}\stackrel{ind.}{\sim}\mathrm{NEF\mbox{-}QVF}[\theta
_{i},V(\theta_{i})/\tau_{i}]$, $i=1,\ldots,p$, be non-Gaussian.
Under the
respective conditions listed in Table~\ref{tab1}, we have
\[
\sup_{0\leq\gamma<\infty}\bigl\llvert \operatorname{URE}^{\mathrm{PG}}(\gamma
)-l_{p}\bigl(\bolds{\theta},\hat{\bolds{\theta}}{}^{\gamma,\bar{Y}}\bigr)
\bigr\rrvert \rightarrow\qquad\mbox{in }L^{1}\mbox{ and in probability, as
}p\rightarrow \infty.
\]
\end{theorem}

\begin{theorem}
\label{theorem8}Let $Y_{i}\stackrel{{ind.}}{\sim}\mathrm
{NEF\mbox{-}QVF}[\theta
_{i},V(\theta_{i})/\tau_{i}]$, $i=1,\ldots,p$, be non-Gaussian. Assume
the respective conditions listed in Table~\ref{tab1}. Then for any estimator
$\hat{\bolds{\theta}}{}^{\hat{\gamma},\bar{Y}}=\frac{\bolds{\tau}}{\bolds
{\tau}+\hat{%
\gamma}}\mathbf{Y}+\frac{\hat{\gamma}}{\bolds{\tau}+\hat{\gamma}}\bar
{Y}$, where $%
\hat{\gamma}\geq0$, we have
\[
\lim_{p\rightarrow\infty}P \bigl( l_{p}\bigl(\bolds{\theta},\hat{
\bolds {\theta}}%
{}^{\mathrm{PG}}\bigr)\geq l_{p}\bigl(\bolds{
\theta},\hat{\bolds{\theta}}{}^{\hat{\gamma},\bar
{Y}%
}\bigr)+\varepsilon \bigr) =0\qquad\mbox{for
any } \varepsilon>0
\]
and
\[
\liminf_{p\rightarrow\infty} \bigl[ R_{p}\bigl(\bolds{\theta},
\hat{%
\bolds{\theta}}{}^{\mathrm{PG}}\bigr)-R_{p}\bigl(\bolds{
\theta},\hat{\bolds{\theta}}{}^{\hat
{\gamma},\bar{Y}%
}\bigr) \bigr] \leq0.
\]
\end{theorem}

\section{Simulation study}\label{sec5}

In this section, we conduct a number of simulations to investigate the
performance of the URE estimators and compare their performance to that of
other existing shrinkage estimators. For each simulation, we first draw
$%
(\theta_{i},\tau_{i})$, $i=1,\ldots,p$, independently from a distribution
and then generate $Y_{i}$ given $(\theta_{i},\tau_{i})$. This process is
repeated a large number of times ($N=100\mbox{,}000$) to obtain an accurate
estimate of the risk for each estimator. The sample size $p$ is chosen to
vary from $20$ to $500$ at an interval of length $20$. For notational
convenience, in this section, we write $A_{i}=1/\tau_{i}$ so that $A_{i}$
is (essentially) the variance.

\subsection{Location-scale family}\label{sec5.1}

For the location-scale families, we consider three non-Gaussian cases: the
Laplace [where the standard variate $Z$ has density $f(z)=\frac
{1}{2}\exp
(-|z|)$], logistic [where the standard variate $Z$ has density $%
f(z)=e^{-z}/(1+e^{-z})^{2}$] and Student-$t$ distributions with 7 degrees
of freedom. To evaluate the performance of the semiparametric URE estimator
$\hat{\theta}{}^{\mathrm{SM}}$, we compare it to the naive estimator
\[
\hat{\bolds{\theta}}{}^{\mathrm{Naive}}_{i}=Y_{i}
\]
and the extended James--Stein estimator
%
\begin{equation}
\hat{\bolds{\theta}}{}^{\mathrm{JS}+}_{i}=\hat{\mu}{}^{\mathrm{JS}+}+
\biggl(1-\frac{(p-3)}{%
\sum_{i=1}^{p}(Y_{i}-\hat{\mu}{}^{\mathrm{JS}+})^{2}/A_{i}} \biggr)^{+}\cdot Y_{i},
\label{extendedJS}
\end{equation}
where $A_{i}=1/\tau_{i}$ and $\hat{\mu}{}^{\mathrm{JS}+}=\sum_{i=1}^{p}(X_{i}/A_{i})/%
\sum_{i=1}^{p}1/A_{i}$.

For each of the three distributions we study four different setups to
generate $(\theta_{i},A_{i}=1/\tau_{i})$ for $i=1,\ldots,p$. We then
generate $Y_{i}$ via (\ref{loc-scale}) except for scenario (4) below.

\begin{figure}

\includegraphics{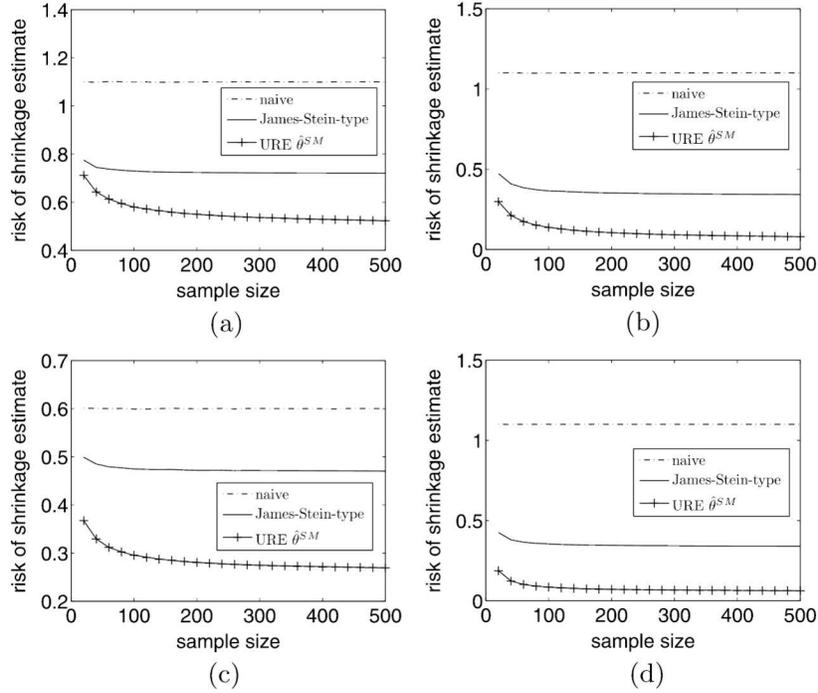}

\caption{Comparison of the risks of different shrinkage estimators for the
Laplace case. \textup{(a)}~$A\sim\operatorname{Unif}(0.1,1)$, $\theta\sim N(0,1)$
independently; $Y=\theta+ \sqrt{A}\cdot Z$. \textup{(b)} $A\sim\operatorname{Unif}(0.1,1)$, $\theta=A$; $Y= \theta+
\sqrt{A}\cdot Z$. \textup{(c)} $A\sim\frac{1}{2}\cdot1_{\{A=0.1\}}+\frac{1}{2}\cdot
1_{\{A=0.5\}}$, $\theta|A=0.1\sim N(2,0.1)$, $\theta|A=0.5\sim
N(0,0.5)$; $Y=\theta+
\sqrt{A}\cdot Z$. \textup{(d)} $A\sim\operatorname{Unif}(0.1,1)$, $\theta=A$; $Y \sim
\operatorname{Unif}[\theta-\sqrt{6A},
\theta+\sqrt{6A}]$.}
\label{laplace}
\end{figure}

\begin{figure}

\includegraphics{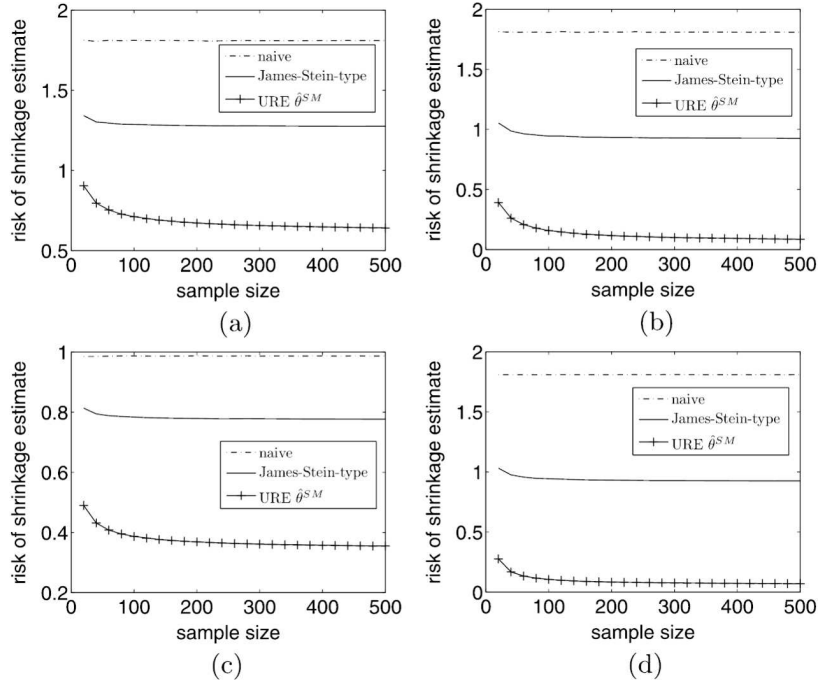}

\caption{Comparison of the risks of different shrinkage estimators for the
logistic case. \textup{(a)}~$A\sim\operatorname{Unif}(0.1,1)$, $\theta\sim N(0,1)$
independently; $Y=\theta+ \sqrt{A}\cdot Z$.
\textup{(b)} $A\sim\operatorname{Unif}(0.1,1)$, $\theta=A$; $Y= \theta+
\sqrt{A}\cdot Z $.
\textup{(c)} $A\sim\frac{1}{2}\cdot1_{\{A=0.1\}}+\frac{1}{2}\cdot
1_{\{A=0.5\}}$, $\theta|A=0.1\sim N(2,0.1)$, $\theta|A=0.5\sim
N(0,0.5)$; $Y=\theta+
\sqrt{A}\cdot Z$.
\textup{(d)} $A\sim\operatorname{Unif}(0.1,1)$, $\theta=A$; $Y \sim
\operatorname{Unif}[\theta-\pi\sqrt{A},
\theta+\pi\sqrt{A}]$.}
\label{logistic}
\end{figure}

\begin{figure}

\includegraphics{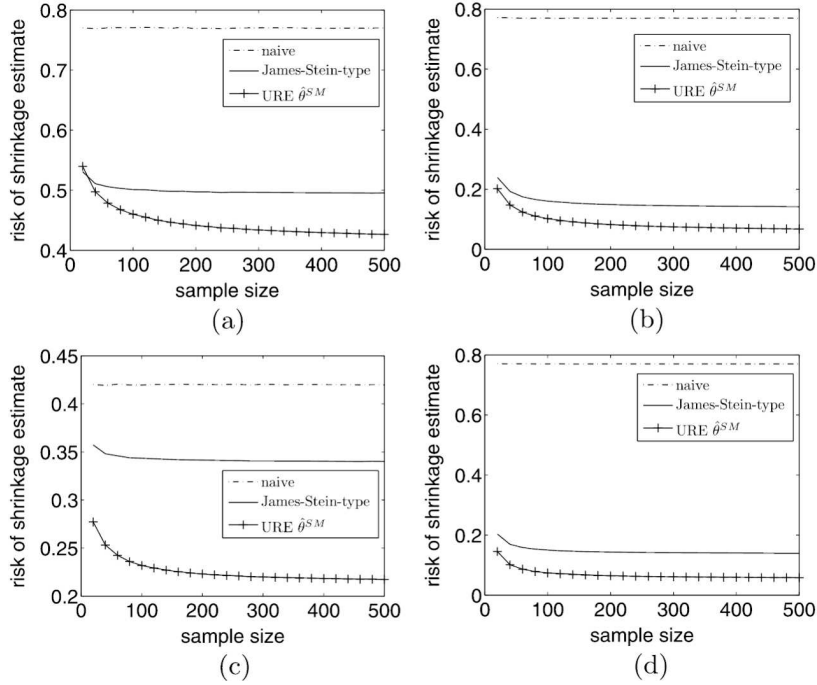}

\caption{Comparison of the risks of different shrinkage estimators for
the $%
t $-distribution ($df = 7$).
\textup{(a)} $A\sim\operatorname{Unif}(0.1,1)$, $\theta\sim N(0,1)$
independently; $Y=\theta+ \sqrt{A}\cdot Z
$.
\textup{(b)} $A\sim\operatorname{Unif}(0.1,1)$, $\theta=A$; $Y= \theta+
\sqrt{A}\cdot Z
$.
\textup{(c)} $A\sim\frac{1}{2}\cdot1_{\{A=0.1\}}+\frac{1}{2}\cdot
1_{\{A=0.5\}}$, $\theta|A=0.1\sim N(2,0.1)$, $\theta|A=0.5\sim
N(0,0.5)$; $Y=\theta+
\sqrt{A}\cdot Z$.
\textup{(d)} $A\sim\operatorname{Unif}(0.1,1)$, $\theta=A$; $Y \sim
\operatorname{Unif}[\theta-\sqrt{21A/5},
\theta+\sqrt{21A/5}]$.}
\label{tdf7}
\end{figure}

\textit{Scenario} (1). $(\theta_{i},A_{i})$ are drawn from $A_{i}\sim%
\operatorname{Unif}(0.1,1)$ and $\theta_{i}\sim N(0,1)$ independently. In this
scenario, the location and scale are independent of each other. Panels (a)
in Figures \ref{laplace}--\ref{tdf7} plot the performance of the three
estimators. The risk
function of the naive estimator $\hat{\bolds{\theta}}{}^{\mathrm{Naive}}_{i}$,
being a
constant for all $p$, is way above the other two. The risk of the
semiparametric URE estimator is significantly smaller than that of
both the
extended James--Stein estimator and the naive estimator, particularly
so when
the sample size $p>40$.

\textit{Scenario} (2). $(\theta_{i},A_{i})$ are drawn from $A_{i}\sim%
\operatorname{Unif}(0.1,1)$ and $\theta_{i}=A_{i}$. This scenario
tests the case
that the location and scale have a strong correlation. Panels (b) in
Figures \ref{laplace}--\ref{tdf7} show the performance of the three
estimators. The risk of the
semiparametric URE estimator is significantly smaller than that of
both the
extended James--Stein estimator and the native estimator. The naive estimator
$\hat{\bolds{\theta}}{}^{\mathrm{Naive}}_{i}$, with a constant risk, performs the worst.
This example indicates that the semiparametric URE estimator behaves
robustly well even when there is a strong correlation between the location
and the scale. This is because the semiparametric URE estimator does not
make any assumption on the relationship between $\theta_{i}$ and $\tau
_{i}$.

\textit{Scenario} (3). $(\theta_{i},A_{i})$ are drawn such that
$A_{i}\sim
\frac{1}{2}\cdot1_{\{A_{i}=0.1\}}+\frac{1}{2}\cdot1_{\{A_{i}=0.5\}}$---that is, $A_{i}$ is 0.1 or 0.5 with 50\% probability each---and that
conditioning on $A_{i}$ being 0.1 or 0.5, $\theta_{i}|A_{i}=0.1\sim
N(2,0.1) $; $\theta_{i}|A_{i}=0.5\sim N(0,0.5)$. This scenario tests the
case that there are two underlying groups in the data. Panels (c) in
Figures~\ref{laplace}--\ref{tdf7} compare the performance of the
semiparametric URE estimator to that of
the native estimator and the extended James--Stein estimator. The
semiparametric URE estimator is seen to significantly outperform the other
two estimators.

%
%

%

\textit{Scenario} (4). $(\theta_{i},A_{i})$ are drawn from $A_{i}\sim%
\operatorname{Unif}(0.1,1)$ and $\theta_{i}=A_{i}$. Given $\theta
_{i}$ and $A_{i}$, $Y_{i}$ are generated from $Y_{i}\sim\operatorname{Unif}(\theta
_{i}-\sqrt{3A_{i}}%
\sigma,\theta+\sqrt{3A_{i}}\sigma)$, where $\sigma$ is the standard
deviation of the standard variate $Z$, that is, $\sigma=\sqrt{2}$,
$\pi/\sqrt{%
3}$ and $\sqrt{7/5}$ for the Laplace, logistics and $t$-distribution
($df =
7$), respectively. This scenario tests the case of model mis-specification
and, hence, the robustness of the estimators. Note that $Y_{i}$ is drawn
from a uniform distribution, not from the Laplace, logistic or $t$
distribution. Panels (d) in Figures \ref{laplace}--\ref{tdf7} show the
performance of the three
estimators. It is seen that the naive estimator behaves the worst and that
the semiparametric URE estimator clearly outperforms the other two. This
example indicates the robust performance of the semiparametric URE
estimator even when the model is incorrectly specified. This is because URE
estimator essentially only involves the first two moments of $Y_{i}$; it
does not rely on the specific density function of the distribution.

%
%
%
%

%



%

\subsection{Exponential family}\label{sec5.2}

We consider exponential family in this subsection, conducting simulation
evaluations on the beta-binomial and Poisson-gamma models.

\subsubsection{Beta-binomial hierarchical model}\label{sec5.2.1}

For binomial observations $Y_{i}\stackrel{\mathrm{ind.}}{\sim
}\operatorname{Bin}(\tau_{i},\theta
_{i})/\tau_{i}$, classical hierarchical inference typically assumes the
conjugate prior $\theta_{i}\stackrel{\mathrm{i.i.d.}}{\sim
}\operatorname{Beta}(\alpha,\beta)$. The
marginal distribution of $Y_{i}$ is used by classical empirical Bayes
methods to estimate the hyper-parameters. Plugging the estimate of these
hyper-parameters into the posterior mean of $\theta_{i}$ given $Y_{i}$
yields the empirical Bayes estimate of $\theta_{i}$. In this
subsection, we
consider both the semiparametric URE estimator
$\hat{\bolds{\theta}}{}^{\mathrm{SM}}$
[equation (\ref{semiURE})] and the parametric URE estimator $%
\hat{\bolds{\theta}}{}^{\mathrm{PM}}$ [equation~(\ref{parURE})], and compare them
with the
parametric empirical Bayes maximum likelihood estimator $\hat{\bolds
{\theta}}%
{}^{\mathrm{ML}}$ and the parametric empirical Bayes method-of-moment estimator $%
\hat{\bolds{\theta}}{}^{\mathrm{MM}}$. The empirical Bayes maximum likelihood
estimator $%
\hat{\bolds{\theta}}{}^{\mathrm{ML}}$ is given by
\[
\hat{\theta}{}^{\mathrm{ML}}_{i}=\frac{\tau_{i}}{\tau_{i}+
\hat{\gamma}{}^{\mathrm{ML}}}\cdot
Y_{i}+\frac{\hat{\gamma}{}^{\mathrm{ML}}}{\tau_{i}+\hat{\gamma}{}^{\mathrm{ML}}}\cdot\hat {\mu}%
{}^{\mathrm{ML}},
\]
where $(\hat{\gamma}{}^{\mathrm{ML}},\hat{\mu}{}^{\mathrm{ML}})$ maximizes the marginal likelihood
of $Y_{i}$:
\[
\bigl(\hat{\gamma}{}^{\mathrm{ML}},\hat{\mu}{}^{\mathrm{ML}}\bigr)=\arg\max
_{\gamma\geq0,\mu
}\prod_{i}\frac{\Gamma(\gamma\mu+\tau_{i}y_{i})\Gamma(\gamma
(1-\mu)+(1-y_{i})\tau_{i})\Gamma(\gamma)}{\Gamma(\gamma+\tau
_{i})\Gamma(\gamma\mu)\Gamma(\gamma(1-\mu))},
\]
where $\mu=\alpha/(\alpha+\beta)$ and $\gamma=\alpha+\beta$ as in
Table~\ref{tab2}. Likewise, the empirical Bayes method-of-moment
estimator $%
\hat{\bolds{\theta}}{}^{\mathrm{MM}}$ is given by
\[
\hat{\theta}{}^{\mathrm{MM}}_{i}=\frac{\tau_{i}}{\tau_{i}+\hat{\gamma}{}^{\mathrm{MM}}}\cdot
Y_{i}+\frac{\hat{\gamma}{}^{\mathrm{MM}}}{\tau_{i}+\hat{\gamma}{}^{\mathrm{MM}}}\cdot\hat {\mu}%
{}^{\mathrm{MM}},
\]
where 
%
\begin{eqnarray*}
\hat{\mu}{}^{\mathrm{MM}} &=&\bar{Y}=\frac{1}{p}\sum
_{i=1}^{p}Y_{i},
\\
\hat{\gamma}{}^{\mathrm{MM}} &=&\frac{\bar{Y}(1-\bar{Y})\cdot
\sum_{i=1}^{p}(1-1/\tau_{i})}{ [ \sum_{i=1}^{p} (
Y_{i}^{2}-\bar{Y}/\tau_{i}-\bar{Y}^{2}(1-1/\tau_{i}) )  ]
^{+}}.
\end{eqnarray*}

There are in total four different simulation setups in which we study the
four different estimators. In addition, in each case, we also calculate the
oracle \emph{risk} ``estimator'' $%
\tilde{\bolds{\theta}}{}^{\mathrm{OR}}$, defined as
%
\begin{equation}
\tilde{\theta}{}^{\mathrm{OR}}_{i}=\frac{\tau_{i}}{\tau_{i}+\tilde{\gamma
}{}^{\mathrm{OR}}}\cdot
Y_{i}+\frac{\tilde{\gamma}{}^{\mathrm{OR}}}{\tau_{i}+\tilde{\gamma}{}^{\mathrm{OR}}}\cdot \tilde{%
\mu}{}^{\mathrm{OR}},
\label{oracle}
\end{equation}
where
\begin{eqnarray*}
\bigl(\tilde{\gamma}{}^{\mathrm{OR}},\tilde{\mu}{}^{\mathrm{OR}}\bigr)&=&\arg\min
_{\gamma\geq0, \mu
}R_{p}\bigl(\bolds{\theta},\hat{\bolds{
\theta}}{}^{\gamma,\mu}\bigr)
\\
&=&\arg\min_{\gamma\geq0, \mu}\sum_{i=1}^{p}
\frac{1}{p}E \biggl[ \biggl( \frac{\tau_{i}}{\tau_{i}+\gamma}Y_{i}+
\frac{\gamma}{\tau
_{i}+\gamma}\mu-\theta_{i} \biggr) ^{2} \biggr].
\end{eqnarray*}
Clearly, the oracle risk estimator $\tilde{\bolds{\theta}}{}^{\mathrm{OR}}$
cannot be used
in practice, since it depends on the unknown $\bolds{\theta}$, but it does
provide a sensible lower bound of the risk achievable by any shrinkage
estimator with the given parametric form.

\begin{example}\label{ex1}
We generate $\tau_{i}\sim\operatorname{
Poisson}(3)+2$ and $\theta_{i}\sim\operatorname{Beta}(1,1)$
independently, and
draw $Y_{i}\sim\operatorname{Bin}(\tau_{i},\theta_{i})/\tau_{i}$.
The oracle
estimator $\tilde{\bolds{\theta}}{}^{\mathrm{OR}}$ is found to have $\tilde{\gamma}{}^{\mathrm{OR}}=2$
and $\tilde{\mu}{}^{\mathrm{OR}}=0.5$. The corresponding risk for the oracle estimator
is numerically found to be $R_{p}(\bolds{\theta
},\tilde{\bolds{\theta}}{}^{\mathrm{OR}})\approx0.0253$. The plot in Figure~\ref%
{riskplots}(a) shows the risks of the five shrinkage estimators as the
sample size $p$ varies. It is seen that the performance of all four
shrinkage estimators approaches that of the parametric oracle
estimator, the
``best estimator'' one can hope to get
under the parametric form. Note that the two empirical Bayes estimators
converges to the oracle estimator faster than the two URE shrinkage
estimators. This is because the hierarchical distribution on $\tau
_{i}$ and
$\theta_{i}$ are exactly the one assumed by the empirical Bayes estimators.
In contrast, the URE estimators do not make any assumption on the
hierarchical distribution but still achieve rather competitive performance.
When the sample size is moderately large, all four estimators well approach
the limit given by the parametric oracle estimator.
\end{example}

\begin{example}\label{ex2}We generate $\tau_{i}\sim\operatorname{
Poisson}(3)+2$ and $\theta_{i}\sim\frac{1}{2}\operatorname
{Beta}(1,3)+\frac{1}{2}%
\operatorname{Beta}(3,1)$ independently, and draw $Y_{i}\sim
\operatorname{Bin}(\tau
_{i},\theta_{i})/\tau_{i}$. In this example, $\theta_{i}$ no longer comes
from a beta distribution, but $\theta_{i}$ and $\tau_{i}$ are still
independent. The oracle estimator is found to have $\tilde{\gamma}%
{}^{\mathrm{OR}}\approx1.5$ and $\tilde{\mu}{}^{\mathrm{OR}}=0.5$. The corresponding risk
for the
oracle estimator $\hat{\bolds{\theta}}{}^{\tau_{0},\theta}$ is $R_{p}(%
\bolds{\theta},\tilde{\bolds{\theta}}{}^{\mathrm{OR}})\approx0.0248$. The plot in
Figure~\ref%
{riskplots}(b) shows the risks of the five shrinkage estimators as the
sample size $p$ varies. Again, as $p$ gets large, the performance of all
shrinkage estimators eventually approaches that of the oracle estimator.
This observation indicates that the parametric form of the prior on
$\theta
_{i}$ is not crucial as long as $\tau_{i}$ and $\theta_{i}$ are
independent.
\end{example}

\begin{figure}[b]

\includegraphics{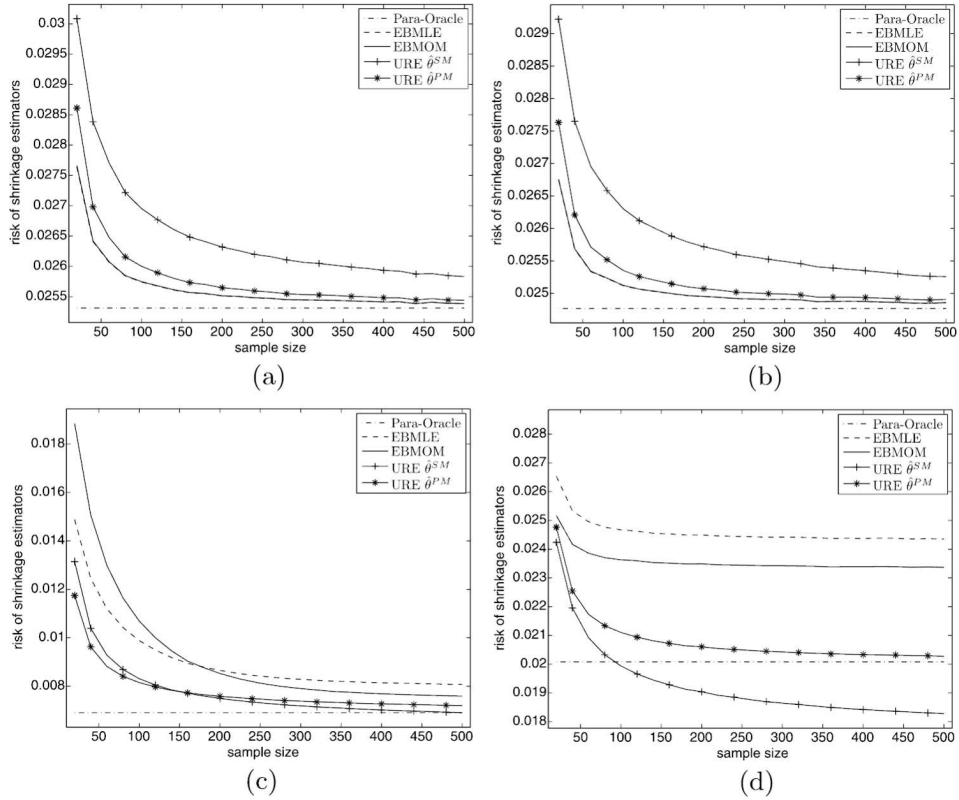}

\caption{Comparison of the risks of shrinkage estimators in the
Beta--Binomial hierarchical models.
\textup{(a)} $\tau\sim\operatorname{Poisson}(3)+2$, $\theta\sim
\operatorname{Beta}(1,1)$ independently;
$Y\sim\operatorname{Bin}(\tau,\theta)/\tau$.
\textup{(b)} $\tau\sim\operatorname{Poisson}(3)+2$, $\theta\sim
\frac{1}{2}\cdot\operatorname{Beta}(1,3)+\frac{1}{2}\cdot
\operatorname{Beta}(3,1)$
independently;
$Y\sim\operatorname{Bin}(\tau,\theta)/\tau$.
\textup{(c)} $\tau\sim\operatorname{Poisson}(3)+2$, $\theta=1/\tau$,
$Y\sim\operatorname{Bin}(\tau,\theta)/\tau$.
\textup{(d)} $I\sim \operatorname{Bern}(1/2)$, $\tau\sim
{I}\cdot\operatorname{Poisson}(10)+{(1-I)}\cdot
\operatorname{Poisson}(1)+2$,
$\theta\sim{I}\cdot\operatorname{Beta}(1,3)+(1-I)\cdot
\operatorname{Beta}(3,1)$;
$Y\sim\operatorname{Bin}(\tau,\theta)/\tau$.}
\label{riskplots}
\end{figure}

\begin{example}\label{ex3}We generate $\tau_{i}\sim\operatorname{
Poisson}(3)+2$ and let $\theta_{i}=1/\tau_{i}$, and then we draw $%
Y_{i}\sim\operatorname{Bin}(\tau_{i},\theta_{i})/\tau_{i}$. In this
case, there
is a (negative) correlation between $\theta_{i}$ and $\tau_{i}$. The
parametric oracle estimator is found to have $\tilde{\gamma}{}^{\mathrm{OR}}\approx
23.0898$ and $\tilde{\mu}{}^{\mathrm{OR}}\approx0.2377$ numerically; the corresponding
risk is $R_{p}(\bolds{\mu
},\tilde{\bolds{\theta}}{}^{\mathrm{OR}})\approx0.0069$. The plot in Figure~\ref%
{riskplots}(c) shows the risks of the five shrinkage estimators as functions
of the sample size $p$. Unlike the previous examples, the two empirical
Bayes estimators no longer converge to the parametric oracle estimator,
that is, the limit of their risk (as $p\rightarrow\infty$) is \emph{strictly}
above the risk of the parametric oracle estimator. On the other hand, the
risk of the parametric URE estimator $\hat{\bolds{\theta}}{}^{\mathrm{PM}}$ still
converges to the risk of the parametric oracle estimator. It is interesting
to note that the limiting risk of the semiparametric URE estimators $%
\hat{\bolds{\theta}}{}^{\mathrm{SM}}$ is actually \emph{strictly} smaller than the
risk of
the parametric oracle estimator (although the difference between the
two is
not easy to spot due to the scale of the plot).
\end{example}

\begin{example}\label{ex4}We generate $(\tau_{i},\theta
_{i},Y_{i})$ as follows. First, we draw $I_{i}$ from $\operatorname{Bernoulli}(1/2)$,
and then
generate $\tau_{i}\sim I_{i}\cdot\operatorname
{Poisson}(10)+(1-I_{i})\cdot\operatorname{%
Poisson}(1)+2$ and $\theta_{i}\sim I_{i}\cdot\operatorname{Beta}%
(1,3)+(1-I_{i})\cdot\operatorname{Beta}(3,1)$. Given $(\tau
_{i},\theta_{i})$, we
draw $Y_{i}\sim\operatorname{Bin}(\tau_{i},\theta_{i})/\tau_{i}$.
In this
example, there exist two groups in the data (indexed by $I_{i}$). It thus
serves to test the different estimators in the presence of grouping. The
parametric oracle estimator is found to have $\tilde{\gamma}{}^{\mathrm{OR}}\approx
0.3108$ and $\tilde{\mu}{}^{\mathrm{OR}}\approx2.0426$; the corresponding risk is
$%
R_{p}(\bolds{\mu},\tilde{\bolds{\theta}}{}^{\mathrm{OR}})\approx0.0201$. Figure~\ref%
{riskplots}(d) plots the risks of the five shrinkage estimators versus the
sample size $p$. The two empirical Bayes estimators clearly encounter much
greater risk than the URE estimators, and the limiting risks of the two
empirical Bayes estimators are significantly larger than the risk of the
parametric oracle estimator. The risk of the parametric URE estimator $%
\hat{\bolds{\theta}}{}^{\mathrm{PM}}$ converges to that of the parametric oracle
estimator. It is quite noteworthy that the semiparametric URE
estimator $%
\hat{\bolds{\theta}}{}^{\mathrm{SM}}$ achieves a significant improvement over the
parametric oracle one.
\end{example}

%
%
%

%
%

\subsubsection{Poisson--Gamma hierarchical model}\label{sec5.2.2}

For Poisson observations\break $Y_{i}\stackrel{\mathrm{ind.}}{\sim
}\operatorname{Poisson}(\tau
_{i}\theta_{i})/\tau_{i}$, the conjugate prior is $\theta
_{i}\stackrel{%
\mathrm{i.i.d.}}{\sim}\Gamma(\alpha,\lambda)$. Like in the previous
subsection, we
compare five estimators: the empirical Bayes maximum likelihood estimator,
the empirical Bayes method-of-moment estimator, the parametric and
semiparametric URE estimators and the parametric oracle ``estimator''
(\ref{oracle}). The empirical Bayes maximum
likelihood estimator $\hat{\bolds{\theta}}{}^{\mathrm{ML}}$ is given by
\[
\hat{\theta}{}^{\mathrm{ML}}_{i}=\frac{\tau_{i}}{\tau_{i}+\hat{\gamma}{}^{\mathrm{ML}}}\cdot
Y_{i}+\frac{\hat{\gamma}{}^{\mathrm{ML}}}{\tau_{i}+\hat{\gamma}{}^{\mathrm{ML}}}\cdot\hat {\mu}%
{}^{\mathrm{ML}},
\]
where $(\hat{\gamma}{}^{\mathrm{ML}},\hat{\mu}{}^{\mathrm{ML}})$ maximizes the marginal likelihood
of $Y_{i}$:
\[
\bigl(\hat{\gamma}{}^{\mathrm{ML}},\hat{\mu}{}^{\mathrm{ML}}\bigr)=\arg\max
_{\gamma\geq0,\mu
}\prod_{i}\frac{\gamma^{\gamma\mu}\Gamma(\gamma\mu+\tau
_{i}y_{i})}{(\tau_{i}+\gamma)^{\tau_{i}y_{i}+\gamma\mu}\Gamma
(\gamma
\mu)},
\]
where $\mu=\alpha\lambda$ and $\gamma=1/\lambda$ as in Table~\ref
{tab2}. The empirical Bayes method-of-moment estimator $\hat{\bolds{\theta
}}{}^{\mathrm{MM}}$ is
given by
\[
\hat{\theta}{}^{\mathrm{MM}}_{i}=\frac{\tau_{i}}{\tau_{i}+\hat{\gamma}{}^{\mathrm{MM}}}\cdot
Y_{i}+\frac{\hat{\gamma}{}^{\mathrm{MM}}}{\tau_{i}+\hat{\gamma}{}^{\mathrm{MM}}}\cdot\hat {\mu}%
{}^{\mathrm{MM}},
\]
where
\begin{eqnarray*}
\hat{\mu}{}^{\mathrm{MM}} &=&\bar{Y}=\frac{1}{p}\sum
_{i=1}^{p}Y_{i},
\\
\hat{\gamma}{}^{\mathrm{MM}} &=&\frac{p\cdot\bar{Y}}{ [ \sum_{i=1}^{p}%
 ( Y_{i}^{2}-\bar{Y}/\tau_{i}-\bar{Y}^{2} )  ] ^{+}}.
\end{eqnarray*}

We consider four different simulation settings.

\begin{figure}

\includegraphics{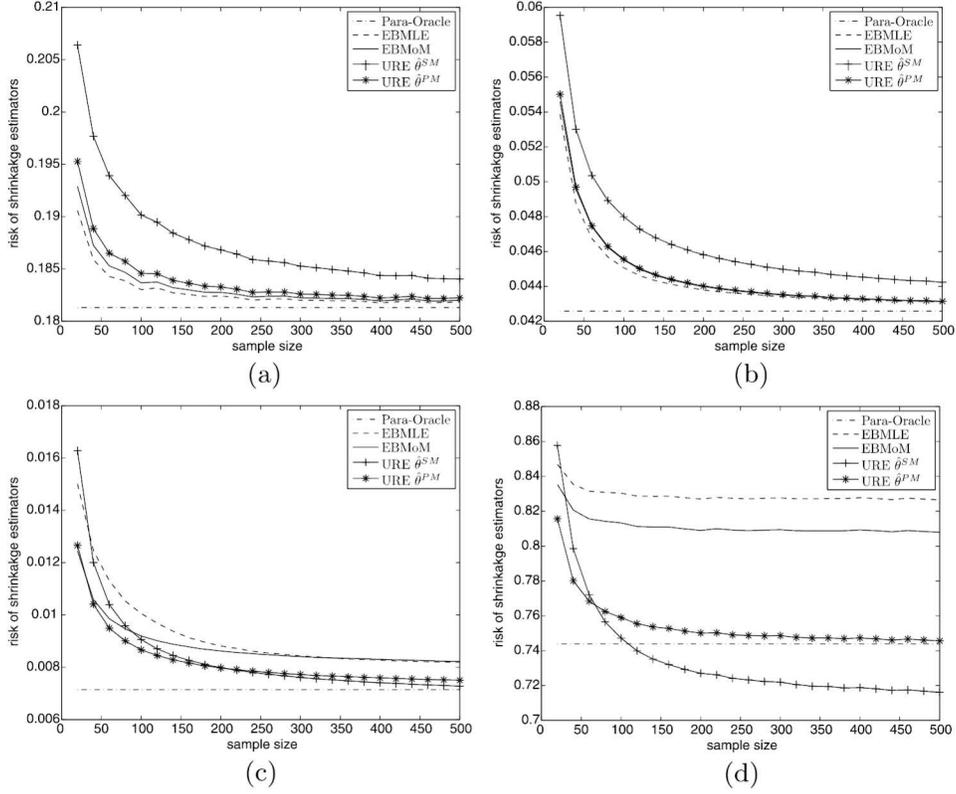}

\caption{Comparison of the risks of shrinkage estimators in Poisson--Gamma
hierarchical model.
\textup{(a)} $\tau\sim\operatorname{Poisson}(3)+2$, $\theta\sim
\Gamma(1,1)$ independently;
$Y\sim\operatorname{Poisson}(\tau\theta)/\tau$.
\textup{(b)}~$\tau\sim\operatorname{Poisson}(3)+2$, $\theta\sim
\operatorname{Unif}(0.1,1)$ independently;
$Y\sim\operatorname{Poisson}(\tau\theta)/\tau$.
\textup{(c)} $\tau\sim\operatorname{Poisson}(3)+2$, $\theta=1/\tau$,
$Y\sim\operatorname{Poisson}(\tau\theta)/\tau$.
\textup{(d)} $I\sim \operatorname{Bern}(1/2)$, $\tau\sim{I}\cdot
\operatorname{Poisson}(10)+(1-I)
\cdot\operatorname{Poisson}(1)+2$, $\theta\sim
{I}\cdot\operatorname{Gamma}(1,1)+(1-I)\cdot
\operatorname{Gamma}(5,1)$;
$Y\sim\operatorname{Poisson}(\tau\theta)/\tau$.}
\label{poissonriskplots}
\end{figure}

\begin{example}\label{ex5}We generate $\tau_{i}\sim\operatorname{
Poisson}(3)+2$ and $\theta_{i}\sim\Gamma(1,1)$ independently, and
draw $%
Y_{i}\sim\operatorname{Poisson}(\tau_{i}\theta_{i})/\tau_{i}$. The
plot in
Figure~\ref{poissonriskplots}(a) shows the risks of the five shrinkage
estimators as the sample size $p$ varies. Clearly, the performance of all
shrinkage estimators approaches that of the parametric oracle
estimator. As
in the beta-binomial case, the two empirical Bayes estimators converge to
the oracle estimator faster than the two URE shrinkage estimators. Again,
this is because the hierarchical distribution on $\tau_{i}$ and $\theta
_{i} $ are exactly the one assumed by the empirical Bayes estimators. The
URE estimators, without making any assumption on the hierarchical
distribution, still achieve rather competitive performance.
\end{example}

\begin{example}\label{ex6} We generate $\tau_{i}\sim\operatorname
{%
Poisson}(3)+2$ and $\theta_{i}\sim\operatorname{Unif}(0.1,1)$
independently, and
draw $Y_{i}\sim\operatorname{Poisson}(\tau_{i}\theta_{i})/\tau
_{i}$. In this
setting, $\theta_{i}$ no longer comes from a gamma distribution, but $%
\theta_{i}$ and $\tau_{i}$ are still independent. The plot in Figure~\ref%
{poissonriskplots}(b) shows the risks of the five shrinkage estimators as
the sample size $p$ varies. As $p$ gets large, the performance of all
shrinkage estimators eventually approaches that of the oracle estimator.
Like in the beta-binomial case, the picture indicates that the parametric
form of the prior on $\theta_{i}$ is not crucial as long as $\tau
_{i}$ and
$\theta_{i}$ are independent.
\end{example}

\begin{example}\label{ex7}We generate $\tau_{i}\sim\operatorname{
Poisson}(3)+2$ and let $\theta_{i}=1/\tau_{i}$, and then we draw $%
Y_{i}\sim\operatorname{Poisson}(\tau_{i}\theta_{i})/\tau_{i}$. In
this setting,
there is a (negative) correlation between $\theta_{i}$ and $\tau
_{i}$. The
plot in Figure~\ref{poissonriskplots}(c) shows the risks of the five
shrinkage estimators as functions of the sample size $p$. Unlike the
previous two examples, the two empirical Bayes estimators no longer converge
to the parametric oracle estimator---the limit of their risk is \emph{%
strictly} above the risk of the parametric oracle estimator. The risk
of the
parametric URE estimator $\hat{\bolds{\theta}}{}^{\mathrm{PM}}$, on the other
hand, still
converges to the risk of the parametric oracle estimator. The limiting risk
of the semiparametric URE estimators $\hat{\bolds{\theta}}{}^{\mathrm{SM}}$ is actually
\emph{strictly} smaller than the risk of the parametric oracle estimator
(although it is not easy to spot it on the plot).
\end{example}

\begin{example}\label{ex8}
We generate $(\tau_{i},\theta_{i})$
by first drawing $I_{i}\sim \operatorname{Bernoulli}(1/2)$ and then $\tau_{i}\sim
I_{i}\cdot\operatorname{Poisson}(10)+(1-I_{i})\cdot\operatorname
{Poisson}(1)+2$ and $%
\theta_{i}\sim I_{i}\cdot\Gamma(1,1)+(1-I_{i})\cdot\Gamma(5,1)$.
With $%
(\tau_{i},\theta_{i})$ obtained, we draw $Y_{i}\sim\operatorname
{Poisson}(\tau
_{i}\theta_{i})/\tau_{i}$. This setting tests the case that there is
grouping in the data. Figure~\ref{poissonriskplots}(d) plots the risks of
the five shrinkage estimators versus the sample size $p$. It is seen that
the two empirical Bayes estimators have the largest risk, and that the
parametric URE estimator $\hat{\bolds{\theta}}{}^{\mathrm{PM}}$ achieves the risk
of the
parametric oracle estimator in the limit. The semiparametric URE estimator
$\hat{\bolds{\theta}}{}^{\mathrm{SM}}$ notably outperforms the parametric oracle
estimator, when $p>100$.
\end{example}

%
%
%
%
%
%

\section{Application to the prediction of batting average}\label{sec6}

In this section, we apply the URE shrinkage estimators to a baseball data
set, collected and discussed in \citet{Bro08}. This data set consists of the
batting records for all the Major League Baseball players in the season of
2005. Following \citet{Bro08}, the data are divided into two half seasons;
the goal is to use the data from the first half season to predict the
players' batting average in the second half season. The prediction can then
be compared against the actual record of the second half season. The
performance of different estimators can thus be directly evaluated.

For each player, let the number of at-bats be $N$ and the successful number
of batting be $H$; we then have
\[
H_{ij}\sim \operatorname{Binomial}(N_{ij},p_{j}),
\]
where $i=1,2$ is the season indicator, $j=1,2,\ldots,$ is the player
indicator, and $p_{j}$ corresponds to the player's hitting ability. Let
$%
Y_{ij}$ be the observed proportion:
\[
Y_{ij}=H_{ij}/N_{ij}.
\]
For this binomial setup, we apply our method to obtain the semiparametric
URE estimators $\hat{\bolds{p}}{}^{\mathrm{SM}}$ and $\hat{\bolds{p}}{}^{\mathrm{SG}}$,
defined in (\ref%
{semiURE}) and (\ref{grandmeanspshrink}), respectively, and the parametric
URE estimators $\hat{\bolds{p}}{}^{\mathrm{PM}}$ and $\hat{\bolds{p}}{}^{\mathrm{PG}}$,
defined in (\ref%
{parURE}) and (\ref{meanUREshrink}), respectively.

To compare the prediction accuracy of different estimators, we note that
most shrinkage estimators in the literature assume normality of the
underlying data. Therefore, for sensible evaluation of different
methods, we
can apply the following variance-stablizing transformation as discussed in
\citet{Bro08}:
\[
X_{ij}=\arcsin\sqrt{\frac{H_{ij}+1/4}{N_{ij}+1/2}},
\]
which gives
\[
X_{ij}\dot{\sim}N\biggl(\theta_{j},\frac{1}{4N_{ij}}
\biggr),\qquad \theta_{j}=\arcsin(%
\sqrt{p_{j}}).
\]
To evaluate an estimator $\hat{\bolds{\theta}}$ based on the
transformed $%
X_{ij} $, we measure the total sum of squared prediction errors (TSE) as
\[
\operatorname{TSE}(\hat{\bolds{\theta}})=\sum_{j}(X_{2j}-
\hat{\theta}%
_{j})^{2}-\sum
_{j}\frac{1}{4N_{2j}}.
\]
To conform to the above transformation (as used by most shrinkage
estimators), we calculate $\hat{\theta}_{j}=\arcsin(\sqrt{\hat{p}_{j}})$,
where $\hat{p}_{j}$ is a URE estimator of the binomial probability, so that
the TSE of the URE estimators can be calculated and compared with other
(normality based) shrinkage estimators.

Table~\ref{tab3} below summarizes the numerical results of our URE
estimators with a collection of competing shrinkage estimators. The values
reported are the ratios of the error of a given estimator to that of the
benchmark naive estimator, which simply uses the first half season $X_{1j}$
to predict the second half $X_{2j}$. All shrinkage estimators are applied
three times---to all the baseball players, the pitchers only, and the
nonpitchers only. The first group of shrinkage estimators in Table~\ref%
{tab3} are the classical ones based on normal theory: two empirical Bayes
methods (applied to $X_{1j}$), the grand mean and the extended James--Stein
estimator (\ref{extendedJS}). The second group includes a number
of more recently developed methods: the
nonparametric shrinkage methods in \citet{BroGre09}, the
binomial mixture model in Muralidharan (\citeyear{Mu10}) and the weighted least squares
and general maximum likelihood estimators (with or without the
covariate of
at-bats effect) in Jiang and Zhang (\citeyear{JiaZha09}, \citeyear{JiaZha10}). The numerical results for
these methods are from \citet{Bro08}, Muralidharan (\citeyear{Mu10}) and Jiang and Zhang
(\citeyear{JiaZha09}, \citeyear{JiaZha10}). The last group corresponds to the results from our binomial
URE methods: the first two are the parametric methods and the last two are
the semiparametric ones.

\begin{figure}[t]

\includegraphics{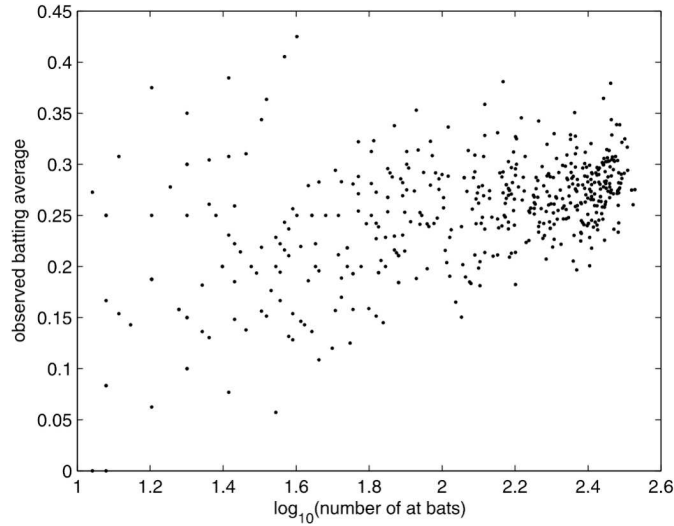}

\caption{Scatter plot of $\log_{10}$ (number of at bats) versus observed
batting average}
\label{figscatter}
\end{figure}

It is seen that our URE shrinkage estimators, especially the semiparametric
ones, achieve very competitive prediction result among all the estimators.
We think the primary reason is that the baseball data contain unique
features that violate the underlying assumptions of the classical empirical
Bayes methods. Both the normal prior assumption and the implicit assumption
of the uncorrelatedness between the binomial probability $p$ and the sample
size $\tau$ are not justified here. To illustrate the last point, we present
a scatter plot of $\log_{10}$ (number of at bats) versus the observed
batting average
$y$ for the nonpitcher group in Figure~\ref{figscatter}.

\begin{table}
\caption{Prediction errors of batting averages for the baseball data}
\label{tab3}
\begin{tabular*}{\textwidth}{@{\extracolsep{\fill}}lccc@{}}
\hline
& \textbf{ALL} & \textbf{Pitchers} & \textbf{Non-Pitchers} \\
\hline
Naive & 1 & 1 & 1 \\[3pt]
Grand mean $\bar{X}_{1\cdot}$ & 0.852 & 0.127 & 0.378 \\
Parametric EB-MM & 0.593 & 0.129 & 0.387 \\
Parametric EB-ML & 0.902 & 0.117 & 0.398 \\
Extended James--Stein & 0.525 & 0.164 & 0.359 \\[3pt]
Nonparametric EB & 0.508 & 0.212 & 0.372 \\
Binomial mixture & 0.588 & 0.156 & 0.314 \\
Weighted least square (Null) & 1.074 & 0.127 & 0.468 \\
Weighted generalized MLE (Null) & 0.306 & 0.173 & 0.326 \\
Weighted least square (AB) & 0.537 & 0.087 & 0.290 \\
Weighted generalized MLE (AB) & 0.301 & 0.141 & 0.261 \\[3pt]
Parametric URE $\hat{\bolds{\theta}}{}^{\mathrm{PG}}$ & 0.515 & 0.105 & 0.278 \\
Parametric URE $\hat{\bolds{\theta}}{}^{\mathrm{PM}}$ & 0.421 & 0.105 & 0.276 \\
Semiparametric URE $\hat{\bolds{\theta}}{}^{\mathrm{SG}}$ & 0.414 & 0.045 & 0.259
\\
Semiparametric URE $\hat{\bolds{\theta}}{}^{\mathrm{SM}}$ & 0.422 & 0.041 & 0.273
\\
\hline
\end{tabular*}
\end{table}

\section{Summary}\label{sec7}

In this paper, we develop a general theory of URE shrinkage estimation in
family of distributions with quadratic variance function. We first
discuss a
class of semiparametric URE estimator and establish their optimality
property. Two specific cases are then carefully studied: the location-scale
family and the natural exponential families with quadratic variance
function. In the latter case, we also study a class of parametric URE
estimators, whose forms are derived from the classical conjugate
hierarchical model. We show that each URE shrinkage estimator is
asymptotically optimal in its own class and their asymptotic optimality do
not depend on the specific distribution assumptions, and more importantly,
do not depend on the implicit assumption that the group mean $\theta$ and
the sample size $\tau$ are uncorrelated, which underlies many classical
shrinkage estimators. The URE estimators are evaluated in comprehensive
simulation studies and one real data set. It is found that the URE
estimators offer numerically superior performance compared to the classical
empirical Bayes and many other competing shrinkage estimators. The
semiparametric URE estimators appear to be particularly competitive.

It is worth emphasizing that the optimality properties of the URE
estimators is not in contradiction with well established results on the
nonexistence of Stein's paradox in simultaneous inference problems with
finite sample space
[\citet{Gut82}], since the results we obtained here are asymptotic ones
when $p$
approaches infinity. A question that naturally arises here is then how
large $p$
needs to be for the URE estimators to become superior compared with
their competitors.
Even though we did not develop a formal finite-sample theory for such
comparison,
our comprehensive simulation indicates that $p$ usually does not need
to be large---$p$ can be as small as 100---for the URE estimators to achieve
competitive performance.

The theory here extends the one on the normal hierarchical models
in \citet{XieKouBro12}. There are three critical features that make
the generalization of previous results possible here: (i) the use of
quadratic risk;
(ii) the linear form of the shrinkage estimator and (iii) the
quadratic variance function of the distribution. The three features together
guarantee the existence of an unbiased risk estimate. For the
hierarchical models
where an unbiased risk estimate does not exist, similar idea can still
be applied
to some estimate of risk, for example, the bootstrap estimate [\citet{Efr04}].
However, a theory is in demand to justify the performance of the
resulting shrinkage estimators.

It would also be an important area of research to study confidence intervals
for the URE estimators obtained here. Understanding whether
the optimality of the URE estimators implies any optimality property of the
estimators of the hyper-parameters under certain conditions is an
interesting related
question. However, such topics are out of scope for the current paper
and we will
need to address them in future research.

\begin{appendix}
\section*{Appendix: Proofs}\label{app}
\setcounter{equation}{0} \renewcommand{\thesection}{A}
\begin{pf*}{Proof of Theorem~\ref{theorem1}} Throughout this proof, when a
supremum is taken, it is over $b_{i}\in{}[0,1]$, $\llvert \mu
\rrvert \leq\max_{i}\llvert  Y_{i}\rrvert $ and Requirement
(MON), unless explicitly stated otherwise. Since
\begin{eqnarray*}
&&\operatorname{URE}(\mathbf{b},\mu)-l\bigl(\bolds{\theta},\hat{\bolds{\theta
}}{}^{\mathbf{b},\mu}\bigr)
\\
&&\qquad=\frac{1}{p}\sum_{i=1}^{p}(1-2b_{i})
\biggl(\frac{V(Y_{i})}{\tau_{i}+\nu
_{2}}-(Y_{i}-\theta_{i})^{2}
\biggr)\\
&&\qquad\quad{}+\frac{2}{p}\sum_{i=1}^{p}b_{i}(Y_{i}-%
\theta_{i}) (\theta_{i}-\mu),
\end{eqnarray*}
it follows that
%
\begin{eqnarray}\label{sure3terms}
\sup\bigl\llvert \operatorname{URE}(\mathbf{b},\mu)-l\bigl(\bolds{\theta},\bolds {
\hat{\theta}}{}^{%
\mathbf{b},\mu}\bigr)\bigr\rrvert &\leq& \frac{1}{p}\Biggl
\llvert \sum_{i=1}^{p}\biggl(%
\frac{V(Y_{i})}{\tau_{i}+\nu_{2}}-(Y_{i}-\theta_{i})^{2}\biggr)
\Biggr\rrvert
\nonumber
\\
&&{}+\frac{2}{p}\sup\Biggl\llvert \sum_{i=1}^{p}b_{i}
\biggl(\frac{V(Y_{i})}{%
\tau_{i}+\nu_{2}}-(Y_{i}-\theta_{i})^{2}
\biggr)\Biggr\rrvert
\\
&&{}+\frac{2}{p}\sup \Biggl\llvert \sum_{i=1}^{p}b_{i}(Y_{i}-
\theta_{i}) (\theta_{i}-\mu )\Biggr\rrvert . \nonumber
\end{eqnarray}
For the first term on the right-hand side, we note that
\begin{eqnarray*}
&&\frac{V(Y_{i})}{\tau_{i}+\nu_{2}}-(Y_{i}-\theta_{i})^{2}\\
&&\qquad=-
\frac{\tau
_{i}%
}{\tau_{i}+\nu_{2}}\bigl(Y_{i}^{2}-EY_{i}^{2}
\bigr)+\biggl(2\theta_{i}-\frac{\nu
_{1}}{%
\tau_{i}+\nu_{2}}\biggr) (Y_{i}-
\theta_{i}).
\end{eqnarray*}
Thus,
\begin{eqnarray*}
&&E \biggl( \biggl(\frac{V(Y_{i})}{\tau_{i}+\nu_{2}}-(Y_{i}-\theta
_{i})^{2}\biggr)^{2} \biggr)
\\
&&\qquad\leq 2 \biggl( \biggl(\frac{\tau_{i}}{\tau_{i}+\nu_{2}}%
\biggr)^{2}
\operatorname{Var}\bigl(Y_{i}^{2}\bigr)+\biggl(2
\theta_{i}-\frac{\nu_{1}}{\tau
_{i}+\nu_{2}%
}\biggr)^{2}\operatorname{Var}(Y_{i})
\biggr)
\\
&&\qquad\leq 2\biggl(\frac{\tau_{i}}{\tau_{i}+\nu_{2}}\biggr)^{2}\operatorname{Var}%
\bigl(Y_{i}^{2}\bigr)+16\theta_{i}^{2}
\operatorname{Var}(Y_{i})+4\biggl(\frac{\nu_{1}}{\tau
_{i}+\nu_{2}}\biggr)^{2}
\operatorname{Var}(Y_{i}).
\end{eqnarray*}
It follows that by conditions (A)--(D)
%
\begin{equation}
\frac{1}{p}\sum_{i=1}^{p} \biggl(
\frac{V(Y_{i})}{\tau_{i}+\nu_{2}}%
-(Y_{i}-\theta_{i})^{2}
\biggr) \rightarrow0\qquad \operatorname{in }L^{2}\mbox{ as }%
p
\rightarrow\infty. \label{varVL2}
\end{equation}
For the term $\frac{2}{p}\sup\llvert \sum_{i}b_{i} ( \frac
{V(Y_{i})}{%
\tau_{i}+\nu_{2}}-(Y_{i}-\theta_{i})^{2} ) \rrvert $ in
(\ref%
{sure3terms}), without loss of generality, let us assume $\tau_{1}\leq
\cdots\leq\tau_{p}$; we then know from Requirement (MON) that
$b_{1}\geq
\cdots\geq b_{2}$. As in Lemma~2.1 in \citet{Li86}, we observe that
\begin{eqnarray*}
&&\sup\frac{2}{p}\Biggl\llvert \sum_{i=1}^{p}b_{i}
\biggl( \frac{V(Y_{i})%
}{\tau_{i}+\nu_{2}}-(Y_{i}-\theta_{i})^{2}
\biggr) \Biggr\rrvert
\\
&&\qquad=\sup_{1\geq b_{1}\geq\cdots\geq b_{p}\geq0}\frac{2}{p}\Biggl\llvert \sum
_{i=1}^{p}b_{i} \biggl(
\frac{V(Y_{i})}{\tau_{i}+\nu_{2}}%
-(Y_{i}-\theta_{i})^{2}
\biggr) \Biggr\rrvert
\\
&&\qquad=\max_{1\leq j\leq p}\frac{2}{p}\Biggl\llvert \sum
_{i=1}^{j} \biggl( \frac{V(Y_{i})}{\tau_{i}+\nu_{2}}%
-(Y_{i}-
\theta_{i})^{2} \biggr) \Biggr\rrvert .
\end{eqnarray*}
Let $M_{j}=\sum_{i=1}^{j}(\frac{V(Y_{i})}{\tau_{i}+\nu
_{2}}-(Y_{i}-\theta
_{i})^{2})$. Then $\{M_{j};j=1,2,\ldots\}$ forms a martingale. The $L^{p}$
maximum inequality implies
\[
E\Bigl(\max_{1\leq j\leq p}M_{j}^{2}\Bigr)\leq 4E
\bigl(M_{p}^{2}\bigr)=4\sum_{i=1}^{p}E
\biggl( \frac{V(Y_{i})}{\tau_{i}+\nu
_{2}}-(Y_{i}-\theta_{i})^{2}
\biggr) ^{2},
\]
which implies by (\ref{varVL2}) that
%
\begin{equation}
\sup\frac{2}{p}\Biggl\llvert \sum_{i=1}^{p}b_{i}
\biggl( \frac{V(Y_{i})}{%
\tau_{i}+\nu_{2}}-(Y_{i}-\theta_{i})^{2}
\biggr) \Biggr\rrvert \rightarrow0%
\qquad\mbox{in }L^{2}\mbox{ as
}p\rightarrow\infty. \label{supVL2}
\end{equation}
For the last term $\frac{2}{p}\sup\llvert \sum_{i}b_{i}(Y_{i}-\theta
_{i})(\theta_{i}-\mu)\rrvert $ in (\ref{sure3terms}), we note that
\[
\frac{1}{p}\sum_{i=1}^{p}b_{i}(Y_{i}-
\theta_{i}) (\theta_{i}-\mu )=\frac{1}{%
p}\sum
_{i=1}^{p}b_{i}
\theta_{i}(Y_{i}-\theta_{i})-
\frac{\mu}{p}%
\sum_{i=1}^{p}b_{i}(Y_{i}-
\theta_{i}).
\]
Using the same argument as in the proof of (\ref{supVL2}), we can show that
\begin{eqnarray*}
\sup\frac{1}{p}\biggl\llvert \sum_{i}b_{i}
\theta_{i}(Y_{i}-\theta _{i})\biggr\rrvert &
\rightarrow&0\qquad\mbox{in }L^{2},
\\
E \biggl( \sup\biggl\llvert \sum_{i}b_{i}(Y_{i}-
\theta_{i})\biggr\rrvert ^{2} \biggr) &=&O(p).
\end{eqnarray*}
Applying condition (E) and the Cauchy--Schwarz inequality, we obtain
\begin{eqnarray*}
&&\frac{1}{p}E \Biggl( \sup\Biggl\llvert \mu \sum
_{i=1}^{p}b_{i}(Y_{i}-
\theta_{i})\Biggr\rrvert \Biggr)
\\
&&\qquad=\frac{1}{p}E \biggl( \max_{1\leq i\leq p}\llvert
Y_{i}\rrvert \cdot \sup\biggl\llvert \sum
_{i}b_{i}(Y_{i}-\theta_{i})
\biggr\rrvert \biggr)
\\
&&\qquad \leq\frac{1}{p} \biggl( E \Bigl( \max_{1\leq i\leq p}Y_{i}^{2}
\Bigr) \cdot E \biggl\{ \sup\biggl\llvert \sum_{i}b_{i}(Y_{i}-
\theta _{i})\biggr\rrvert ^{2} \biggr\} \biggr)
^{1/2}
\\
&&\qquad= O\bigl(p^{(1-\varepsilon+1)/2-1}\bigr)=O\bigl(p^{-\varepsilon/2}\bigr).
\end{eqnarray*}
Therefore,
\[
\frac{1}{p}\sup\Biggl\llvert \mu\sum_{i=1}^{p}b_{i}(Y_{i}-
\theta _{i})\Biggr\rrvert \rightarrow0\mbox{ in }L^{1}.
\]
This completes the proof, since each term on the right-hand side of
(\ref%
{sure3terms}) converges to zero in $L^{1}$.
\end{pf*}

\begin{pf*}{Proof of Theorem~\ref{theorem2}} Throughout this proof, when a
supremum is taken, it is over $b_{i}\in{}[0,1]$, $\llvert \mu
\rrvert \leq\max_{i}\llvert  Y_{i}\rrvert $ and Requirement
(MON). Note that
\[
\operatorname{URE}\bigl(\hat{\mathbf{b}}{}^{{\mathrm{SM}}},\hat{\mu}{}^{{\mathrm{SM}}}\bigr)
\leq\operatorname {URE}(\hat{%
\mathbf{b}},\hat{\mu})
\]
and we know from Theorem~\ref{theorem1} that
\[
\sup\bigl\llvert \operatorname{URE}(\mathbf{b},\mu)-l_{p}\bigl(\bolds{
\theta},\hat {\bolds{\theta}}%
{}^{\mathbf{b},\mu}\bigr)\bigr\rrvert
\rightarrow0\qquad\mbox{in probability.}
\]
It follows that for any $\varepsilon>0$
\begin{eqnarray*}
&&P \bigl( l_{p}\bigl(\bolds{\theta},\hat{\bolds{\theta}}{}^{\mathrm{SM}}
\bigr)\geq l_{p}\bigl(\bolds{\theta},%
\hat{\bolds{
\theta}}{}^{\hat{\mathbf{b}},\hat{\mu}}\bigr)+\varepsilon \bigr)
\\
&&\qquad\leq P \bigl( l_{p}\bigl(\bolds{\theta},\hat{\bolds{
\theta}}{}^{\mathrm{SM}}\bigr)-\operatorname {URE}\bigl(\hat{%
\mathbf{b}}{}^{{\mathrm{SM}}},\hat{\mu}{}^{{\mathrm{SM}}}\bigr)\geq l_{p}
\bigl(\bolds{\theta},\hat {\bolds{\theta}}{}^{%
\hat{\mathbf{b}},\hat{\mu}}\bigr)-\operatorname{URE}(\hat{
\mathbf{b}},\hat{\mu })+\varepsilon \bigr)
\\
&&\qquad\leq P \biggl( \bigl\llvert l_{p}\bigl(\bolds{\theta},\hat{\bolds{
\theta }}{}^{\mathrm{SM}}\bigr)-\operatorname{%
URE}\bigl(\hat{
\mathbf{b}}{}^{{\mathrm{SM}}},\hat{\mu}{}^{{\mathrm{SM}}}\bigr)\bigr\rrvert \geq
\frac
{\varepsilon
}{2} \biggr)
\\
&&\qquad\quad{}+P \biggl( \bigl\llvert l_{p}\bigl(\bolds{\theta},\hat{\bolds{
\theta}}{}^{\hat
{\mathbf{b}},\hat{%
\mu}}\bigr)-\operatorname{URE}(\hat{\mathbf{b}},\hat{\mu})\bigr\rrvert
\geq\frac{%
\varepsilon}{2} \biggr) \rightarrow0.
\end{eqnarray*}
Next, to show that
\[
\limsup_{p\rightarrow\infty} \bigl[ R_{p}\bigl(\bolds{\theta},
\hat{%
\bolds{\theta}}{}^{\mathrm{SM}}\bigr)-R_{p}\bigl(\bolds{
\theta},\hat{\bolds{\theta}}{}^{\hat
{\mathbf{b}},\hat{\mu%
}}\bigr) \bigr] \leq0,
\]
we note that
\begin{eqnarray*}
&&l_{p}\bigl(\bolds{\theta},\hat{\bolds{\theta}}{}^{\mathrm{SM}}
\bigr)-l_{p}\bigl(\bolds{\theta },\hat{%
\bolds{
\theta}}{}^{\hat{\mathbf{b}},\hat{\mu}}\bigr)
\\[-1pt]
&&\qquad= \bigl( l_{p}\bigl(\bolds{\theta},\hat{\bolds{\theta}}{}^{\mathrm{SM}}\bigr)-
\operatorname {URE}\bigl(\hat{\mathbf{b}}%
{}^{{\mathrm{SM}}},\hat{
\mu}{}^{{\mathrm{SM}}}\bigr) \bigr) + \bigl( \operatorname{URE}\bigl(\hat{\mathbf
{b}}{}^{{\mathrm{SM}}},%
\hat{\mu}{}^{{\mathrm{SM}}}\bigr)-\operatorname{URE}(\hat{
\mathbf{b}},\hat{\mu}) \bigr)
\\[-1pt]
&&\qquad\quad{}+ \bigl(\operatorname{URE}(\hat{\mathbf{b}},\hat{\mu})-l_{p}\bigl(\bolds{
\theta },\hat{\bolds{\theta}}{}^{%
\hat{\mathbf{b}},\hat{\mu}}\bigr) \bigr)
\\[-1pt]
&&\qquad\leq 2\sup\bigl\llvert \operatorname{URE}(\mathbf{b},\mu)-l_{p}\bigl(\bolds
{\theta},\hat{%
\bolds{\theta}}{}^{\mathbf{b},\mu}\bigr)\bigr\rrvert .
\end{eqnarray*}
Theorem~\ref{theorem1} then implies that
\[
\limsup_{p\rightarrow\infty} \bigl[ R\bigl(\bolds{\theta},\hat{\bolds {
\theta}}%
{}^{\mathrm{SM}}\bigr)-R\bigl(\bolds{\theta},\hat{\bolds{
\theta}}{}^{\hat{\mathbf{b}},\hat{\mu
}}\bigr) \bigr] \leq0.
\]
\upqed
\end{pf*}

\begin{pf*}{Proof of Theorem~\ref{theorem3}} Throughout this proof, when a
supremum is taken, it is over $b_{i}\in{}[0,1]$ and Requirement (MON).
Since
\begin{eqnarray*}
&&\operatorname{URE}^{G}(\mathbf{b})-l_{p}\bigl(\bolds{\theta},
\hat{\bolds{\theta }}{}^{\mathbf{b},{\bar{%
Y}}}\bigr)
\\[-1pt]
&&\qquad=\frac{1}{p}\sum_{i=1}^{p}
\biggl(1-2\biggl(1-\frac{1}{p}\biggr)b_{i}\biggr) \biggl(
\frac{%
V(Y_{i})}{\tau_{i}+\nu_{2}}-(Y_{i}-\theta_{i})^{2} \biggr)
\\[-1pt]
&&\qquad\quad{}+\frac{2}{p}\sum_{i=1}^{p}b_{i}
\biggl( \theta_{i}(Y_{i}-\theta _{i})+
\frac{1}{p}(Y_{i}-\theta_{i})^{2}-(Y_{i}-
\theta_{i})\bar {Y} \biggr) ,
\end{eqnarray*}
it follows that
%
\begin{eqnarray}\label{ure5terms}
&&\sup\bigl\llvert \operatorname{URE}^{G}(\mathbf{b})-l_{p}\bigl(
\bolds{\theta},%
\hat{\bolds{\theta}}{}^{\mathbf{b},{\bar{Y}}}\bigr)\bigr\rrvert
\nonumber
\\[-1pt]
&&\qquad\leq \frac{1}{p}\biggl\llvert \sum_{i}
\biggl( \frac{V(Y_{i})}{\tau_{i}+\nu_{2}}-(Y_{i}-\theta _{i})^{2}
\biggr) \biggr\rrvert
\nonumber
\\[-1pt]
&&\qquad\quad{}+\frac{2}{p}\biggl(1-\frac{1}{p}\biggr)\sup\biggl\llvert \sum
_{i}b_{i} \biggl( \frac{V(Y_{i})}{\tau_{i}+\nu_{2}}-(Y_{i}-
\theta_{i})^{2} \biggr) \biggr\rrvert
\\[-1pt]
&&\qquad\quad{}+\frac{2}{p}\sup\biggl\llvert \sum_{i}b_{i}
\theta _{i}(Y_{i}-\theta_{i})\biggr\rrvert
\nonumber
\\[-1pt]
&&\qquad\quad{}+\frac{2}{p^{2}}\sum_{i}(Y_{i}-
\theta_{i})^{2}+\frac{2}{p}%
\llvert \bar{Y}
\rrvert \cdot\sup\biggl\llvert \sum_{i}b_{i}(Y_{i}-
\theta_{i})\biggr\rrvert .\nonumber
\end{eqnarray}

We have already shown in the proof of Theorem~\ref{theorem1} that the first
three terms on the right-hand side converge to zero in $L^{2}$. It only
remains to manage the last two terms:
\[
E \biggl( \frac{2}{p^{2}}\sum_{i}(Y_{i}-
\theta_{i})^{2} \biggr) =%
\frac{2}{p^{2}}\sum
_{i}\operatorname{Var}(Y_{i})\rightarrow0
\]
by regularity condition (A),
\[
\frac{1}{p}E \biggl( \llvert \bar{Y}\rrvert \cdot\sup\biggl\llvert \sum
_{i}b_{i}(Y_{i}-
\theta_{i})\biggr\rrvert \biggr) \leq\frac{1}{p}%
E
\biggl( \max_{1\leq i\leq p}\llvert Y_{i}\rrvert \cdot\sup
\biggl\llvert \sum_{i}b_{i}(Y_{i}-
\theta_{i})\biggr\rrvert \biggr) \rightarrow0,
\]
as was shown in the proof of Theorem~\ref{theorem1}. Therefore, the
last two
terms of (\ref{ure5terms}) converge to zero in $L^{1}$, and this completes
the proof.
\end{pf*}

\begin{pf*}{Proof of Theorem~\ref{theorem4}} With Theorem~\ref{theorem3}
established, the proof is almost identical to that of Theorem~\ref
{theorem2}.
\end{pf*}

\begin{pf*}{Proof of Lemma~\ref{theoremLS}} It is straightforward to check that
(i)--(iii) imply conditions (A)--(D) in Section~\ref{sec2}, so we only need
to verify condition (E). Since $Y_{i}^{2}=Z_{i}^{2}/\tau_{i}+\theta
_{i}^{2}+2\theta_{i}Z_{i}/\sqrt{\tau_{i}}$, we know that
%
\begin{equation}
\max_{1\leq i\leq p}Y_{i}^{2}\leq\max
_{i}\frac{1}{\tau_{i}}%
\cdot\max
_{i}Z_{i}^{2}+\max_{i}
\theta _{i}^{2}+2\max_{i}\llvert
\theta_{i}/\sqrt{\tau_{i}}\rrvert \cdot\max
_{i}\llvert Z_{i}\rrvert . \label{maxy2}
\end{equation}
(i) and (ii) imply $\max_{i}1/\tau_{i}=O(p^{1/2})$ and $\max_{i}%
\llvert \theta_{i}/\sqrt{\tau_{i}}\rrvert =O(p^{1/2})$. (iii) gives
$\max_{i}\theta_{i}^{2}=O(p^{2/(2+\varepsilon)})$. We next bound $%
E(\max_{1\leq i\leq p}\llvert  Z_{i}\rrvert ^{k})$ for $k=1,2$.
(iv) implies that for $k=1,2$,
\begin{eqnarray*}
&&E\Bigl(\max_{i}\llvert Z_{i}\rrvert
^{k}\Bigr)
\\[-2pt]
&&\qquad=\int_{0}^{\infty
}kt^{k-1}P\Bigl(\max
_{i}\llvert Z_{i}\rrvert >t\Bigr)\,dt
\\[-2pt]
&&\qquad\leq\int_{0}^{\infty}kt^{k-1}\bigl(1-
\bigl(1-Dt^{-\alpha}\bigr)^{p}\bigr)\,dt
\\[-2pt]
&&\qquad=\int_{0}^{p^{1/\alpha}}kt^{k-1}\bigl(1-
\bigl(1-Dt^{-\alpha
}\bigr)^{p}\bigr)\,dt+\int_{p^{1/\alpha}}^{\infty}kt^{k-1}
\bigl(1-\bigl(1-Dt^{-\alpha}\bigr)^{p}\bigr)\,dt
\\[-2pt]
&&\qquad=O\bigl(p^{k/\alpha}\bigr)+p^{k/\alpha}\int_{1}^{\infty}kz^{k-1}
\biggl( 1-\biggl(1-\frac{%
1}{p}Dz^{-\alpha}\biggr)^{p}
\biggr) \,dz,
\end{eqnarray*}
where a change of variable $z=t/p^{1/\alpha}$ is applied. We know by the
monotone convergence theorem that for $k\geq1$, as $p\rightarrow
\infty$,
\[
\int_{1}^{\infty}z^{k-1} \biggl( 1-\biggl(1-
\frac{1}{p}Dz^{-\alpha}\biggr)^{p} \biggr) \,dz\rightarrow\int
_{1}^{\infty}z^{k-1}\bigl(1-\exp
\bigl(-Dz^{-\alpha
}\bigr)\bigr)\,dz<\infty.
\]
It then follows that
\[
E\Bigl(\max_{1\leq i\leq p}\llvert Z_{i}\rrvert
^{k}\Bigr)=O\bigl(p^{k/\alpha
}\bigr),\qquad\mbox{for }k=1,2.
\]
Taking it back to (\ref{maxy2}) gives
\[
E\Bigl(\max_{1\leq i\leq p}Y_{i}^{2}\Bigr)=O
\bigl(p^{1/2+2/\alpha
}\bigr)+O\bigl(p^{2/(2+\varepsilon
)}\bigr)+O\bigl(p^{1/2+1/\alpha}
\bigr),
\]
which verifies condition (E).\vadjust{\goodbreak}
\end{pf*}

To prove Lemma \ref{lemmaNEFQVF}, we need the following lemma.

\begin{lemma}
\label{lemma1} Let $Y_{i}$ be independent from one of the six NEF-QVFs. Then
condition \textup{(B)} in Section~\ref{sec2} and
\begin{longlist}[(F)]
\item[(F)] $\limsup_{p\rightarrow\infty}\frac{1}{p}%
\sum_{i=1}^{p}\llvert {\theta_i}\rrvert ^{2+\varepsilon
}<\infty$ for some $\varepsilon>0$;

\item[(G)] $\limsup_{p\rightarrow\infty}\frac{1}{p}\sum_{i=1}^{p}%
\operatorname{Var}^{2}(Y_{i})<\infty$;

\item[(H)] $\sup_{i}\operatorname{skew}(Y_{i})=\sup_{i}\frac{1}{\sqrt{\tau
_{i}}}\frac{\nu_{1}+2\nu_{2}\theta_i}{(\nu_{0}+\nu_{1}\theta_i+\nu
_{2}\theta_i^{2})^{1/2}}<\infty$;
\end{longlist}
imply condition (E).
\end{lemma}

\begin{pf*}{Proof of Lemma~\ref{lemma1}} Let us denote $\sigma
_{i}^{2}=\operatorname{Var}(Y_{i})$. We can write $Y_{i}=\sigma_{i}Z_{i}+\theta_{i}$,
where $Z_{i}$ are independent with mean zero and variance one. It follows
from $Y_{i}^{2}=\sigma_{i}^{2}Z_{i}^{2}+\theta_{i}^{2}+2\sigma
_{i}\theta
_{i}Z_{i}$ that
%
\begin{equation}
\max_{1\leq i\leq p}Y_{i}^{2}\leq\max
_{i}\sigma_{i}^{2}\cdot \max
_{i}Z_{i}^{2}+\max_{i}
\theta _{i}^{2}+2\max_{i}
\sigma_{i}\llvert \theta_{i}\rrvert \cdot \max
_{i}\llvert Z_{i}\rrvert . \label{bigOmax}
\end{equation}
Condition (B) implies $\max_{i}\sigma_{i}^{2}\theta_{i}^{2}\leq
\sum_{i}\sigma_{i}^{2}\theta_{i}^{2}=O(p)$. Thus, $\max_{i}\sigma
_{i}\llvert \theta_{i}\rrvert = O(p^{1/2})$. Similarly,
condition (G)
implies that $\max_{i}\sigma_{i}^{2}=O(p^{1/2})$. Condition~(F) implies
that $\max_{i}\llvert {\theta_{i}}\rrvert ^{2+\varepsilon}\leq
\sum_{i}\llvert {\theta_{i}}\rrvert ^{2+\varepsilon}=O(p)$, which
gives $\max_{i}\theta_{i}^{2}=O(p^{2/(2+\varepsilon)})$. If we can show
that
%
\begin{equation}
E\Bigl(\max_{1\leq i\leq p}\llvert Z_{i}\rrvert \Bigr)=O(\log
p),\qquad E\Bigl(\max_{1\leq i\leq p}Z_{i}^{2}\Bigr)=O
\bigl(\log^{2}p\bigr), \label{maxmean}
\end{equation}
then we establish (E), since
\[
E\Bigl(\max_{1\leq i\leq p}Y_{i}^{2}\Bigr)=O
\bigl(p^{1/2}\log^{2}p+p^{2/(2+\varepsilon
)}+p^{1/2}\log p
\bigr)=O\bigl(p^{2/(2+\varepsilon^{\ast})}\bigr),
\]
where $\varepsilon^{\ast}=\min(\varepsilon,1)$. To prove (\ref
{maxmean}), we begin from
%
\begin{equation}
E\Bigl(\max_{i}\llvert Z_{i}\rrvert
^{k}\Bigr)=k\int_{0}^{\infty
}t^{k-1}P
\Bigl(\max_{i}\llvert Z_{i}\rrvert >t\Bigr)\,dt\qquad
\mbox{for all }k>0. \label{momentid}
\end{equation}
The large deviation results for NEF-QVF in Morris [(\citeyear{Mor82}), Section~9] and
condition (H) (i.e., $Y_{i}$ have bounded skewness) imply that for all
$t>1$%
, the tail probabilities $P(\llvert  Z_{i}\rrvert >t)$ are uniformly
bounded exponentially: there exists a constant $c_{0}>0$ such that
\[
P\bigl(\llvert Z_{i}\rrvert >t\bigr)\leq e^{-c_{0}t}\qquad\mbox{for
all }i.
\]
Taking it into (\ref{momentid}), we have
%
\begin{eqnarray}\label{partint}
E\Bigl(\max_{i}\llvert Z_{i}\rrvert
^{k}\Bigr) &\leq&\int_{0}^{\infty
}kt^{k-1}
\bigl(1-\bigl(1-e^{-c_{0}t}\bigr)^{p}\bigr)\,dt
\nonumber
\\[-2pt]
&=&\int_{0}^{\log p/c_{0}}kt^{k-1}\bigl(1-
\bigl(1-e^{-c_{0}t}\bigr)^{p}\bigr)\,dt
\nonumber
\\[-9pt]
\\[-9pt]
\nonumber
&&{}+\int_{\log p/c_{0}}^{\infty}kt^{k-1}\bigl(1-
\bigl(1-e^{-c_{0}t}\bigr)^{p}\bigr)\,dt
\nonumber
\\[-2pt]
&=&O\bigl(\log^{k}p\bigr)+\int_{0}^{\infty}k
\biggl(z+\frac{1}{c_{0}}\log p\biggr)^{k-1} \biggl( 1-\biggl(1-
\frac{1}{p}e^{-c_{0}z}\biggr)^{p} \biggr) \,dz,\nonumber
\end{eqnarray}
where in the last line a change of variable $z=t-\log p/c_{0}$ is applied.
We know by the monotone convergence theorem that for $k\geq1$, as
$p\rightarrow\infty$,
\[
\int_{0}^{\infty}z^{k-1} \biggl( 1-\biggl(1-
\frac{1}{p}e^{-c_{0}z}\biggr)^{p} \biggr) \,dz\rightarrow\int
_{0}^{\infty}z^{k-1}\bigl(1-\exp
\bigl(-e^{-c_{0}z}\bigr)\bigr)\,dz<\infty.
\]
It then follows from (\ref{partint}) that
\[
E\Bigl(\max_{1\leq i\leq p}\llvert Z_{i}\rrvert
^{k}\Bigr)=O\bigl(\log^{k}p\bigr),\qquad
\mbox{for }k=1,2,
\]
which completes our proof.
\end{pf*}

\begin{pf*}{Proof of Lemma~\ref{lemmaNEFQVF}} We go over the five
distributions one by one.

(1) Binomial. Since $\theta_{i}=p_{i}$, $\operatorname{Var}%
(Y_{i})=p_{i}(1-p_{i})/n_{i}$, and $\operatorname{Var}(Y_{i}^{2})\leq
EY_{i}^{4}\leq1$, it is straightforward to verify that $n_{i}\geq2$ for
all $i$ guarantees conditions (A)--(E) in Section~\ref{sec2}.

(2) Poisson. $\operatorname{Var}(Y_{i})=\theta_{i}/\tau_{i}$, and $\operatorname
{Var}%
(Y_{i}^{2})=(4\tau_{i}^{2}\theta_{i}^{3}+6\tau_{i}\theta
_{i}^{2}+\theta
_{i})/\tau_{i}^{3}$. It is straightforward to verify that $\inf_{i}\tau
_{i}>0$, $\inf_{i}\tau_{i}\theta_{i}>0$ and $\sum_{i}\theta_{i}^{3}=O(p)$
imply conditions (A)--(D) in Section~\ref{sec2} and conditions (F)--(H) in
Lemma~\ref{lemma1}.

(3) Negative-binomial. $\theta_{i}=\frac{p_{i}}{1-p_{i}}$, $\operatorname
{Var}%
(Y_{i})=\frac{1}{n_{i}}\frac{p_{i}}{(1-p_{i})^{2}}=\frac
{1}{n_{i}}(\theta
_{i}+\theta_{i}^{2})$, so\break $v_{0}=0$, $\nu_{1}=\nu_{2}=1$. $\operatorname
{Var}%
(Y_{i}^{2})=\frac{1}{n_{i}^{3}(1-p_{i})^{4}}%
(p_{i}+4p_{i}^{2}+6n_{i}p_{i}^{2}+p_{i}^{3}+4n_{i}p_{i}^{3}+4n_{i}^{2}p_{i}^{3})
$. From these, we know that $\sum_{i=1}^{p}\operatorname{Var}^{2}(Y_{i})=\sum_{i}%
\frac{1}{(n_{i}p_{i})^{2}}(\frac{p_{i}}{1-p_{i}})^{4}=\break O(\sum_{i}(\frac
{p_{i}%
}{1-p_{i}})^{4})=O(p)$, which verifies conditions (A) and (G). $%
\sum_{i=1}^{p}\operatorname{Var}(Y_{i})\theta_{i}^{2}=\sum_{i}\frac
{1}{n_{i}p_{i}}%
(\frac{p_{i}}{1-p_{i}})^{4}=O(\sum_{i}(\frac{p_{i}}{1-p_{i}})^{4})=O(p)$,
which verifies condition (B). For condition (C), since $n_{i}\geq1$
and $%
0\leq p_{i}\leq1$, we only need to verify that $\sum_{i}\frac{1}{%
n_{i}^{3}(1-p_{i})^{4}}(p_{i}+6n_{i}p_{i}^{2})=O(p)$. This is true,
since $%
\sum_{i}\frac{1}{n_{i}^{3}(1-p_{i})^{4}}(p_{i}+6n_{i}p_{i}^{2})=\sum_{i}(%
\frac{1}{(n_{i}p_{i})^{3}}(\frac{p_{i}}{1-p_{i}})^{4}+ \frac{6}{%
(n_{i}p_{i})^{2}}(\frac{p_{i}}{1-p_{i}})^{4})=O(\sum_{i}(\frac
{p_{i}}{1-p_{i}%
})^{4})=O(p)$. Condition (D) is automatically satisfied. For condition (F),
consider $\sum_{i=1}^{p}\theta_{i}^{4}$. It is $\sum_{i}(\frac{p_{i}}{%
1-p_{i}})^{4}=O(p)$. For condition~(H), note that $\operatorname{skew}(Y_{i})=\frac{1}{
\sqrt{n_{i}}}\frac{p_{i}+1}{\sqrt{p_{i}}}\leq1+1/\sqrt{n_{i}p_{i}}$.
Thus, $%
\sup_{i}\operatorname{skew}(Y_{i})<\infty$ by~(i).

(4) Gamma. $\theta_{i}=\alpha\lambda_{i} $, $\operatorname
{Var}(Y_{i})=\alpha
\lambda_{i}^{2}/\tau_{i}$, so $v_{0}=\nu_{1}=0$, $\nu_{2}=1/\alpha
$. $%
\operatorname{skew}(Y_{i})=2/\sqrt{\tau_{i}\alpha}$. $\operatorname{Var}(Y_{i}^{2})=\frac{%
\alpha}{\tau_{i}}\lambda_{i}^{4}(\alpha+\frac{1}{\tau
_{i}})(4\alpha+%
\frac{6}{\tau_{i}})$. It is straightforward to verify that (i) and (ii)
imply conditions (A)--(D) in Section~\ref{sec2} and conditions (F)--(H) in
Lemma~\ref{lemma1}.

(5) GHS. $\theta_{i}=\alpha\lambda_{i} $, $\operatorname{Var}(Y_{i})=\alpha
(1+\lambda_{i}^{2})/\tau_{i}$, so $v_{0}=\alpha$, $\nu_{1}=0$, $\nu
_{2}=1/\alpha$. $\operatorname{skew}(Y_{i})=(2/\sqrt{\tau_{i}\alpha})\lambda
_{i}/(1+\lambda_{i}^{2})^{1/2}\leq2/\sqrt{\tau_{i}\alpha}$.\vadjust{\goodbreak} $\operatorname
{Var}%
(Y_{i}^{2})=\frac{2\alpha}{\tau_{i}}(1+\lambda_{i}^{2})(\alpha+\frac
{1}{%
\tau_{i}})\*(\frac{1}{\tau_{i}}+\frac{3}{\tau_{i}}\lambda
_{i}^{2}+2\alpha
\lambda_{i}^{2})$. It is then straightforward to verify that (i) and (ii)
imply conditions (A)--(D) in Section~\ref{sec2} and conditions (F)--(H) in
Lemma~\ref{lemma1}.\vspace*{-2pt}
\end{pf*}

\begin{pf*}{Proof of Theorem~\ref{theorem5}} We note that the set over which
the supremum is taken is a subset of that of Theorem~\ref{theorem1}. The
desired result thus automatically holds.\vspace*{-2pt}
\end{pf*}

\begin{pf*}{Proof of Theorem~\ref{theorem6}} With Theorem~\ref{theorem5}
established, the proof is almost identical to that of Theorem~\ref
{theorem2}.\vspace*{-2pt}
\end{pf*}

\begin{pf*}{Proof of Theorem~\ref{theorem7}} We note that the set over which
the supremum is taken is a subset of that of Theorem~\ref{theorem3}. The
desired result thus holds.\vspace*{-2pt}
\end{pf*}

\begin{pf*}{Proof of Theorem~\ref{theorem8}} The proof essentially follows the
same steps in that of Theorem~\ref{theorem2}.\vspace*{-2pt}
\end{pf*}
\end{appendix}

\section*{Acknowledgments}
The views expressed herein are the authors alone
and are not necessarily the views of Two Sigma Investments, LLC or any
of its affliates.

%

\printaddresses
\end{document}